\documentclass[12pt]{article}
\usepackage[utf8]{inputenc}
\usepackage{amssymb,amsmath,amsfonts,eucal,mathrsfs,amsthm} 
\usepackage[colorinlistoftodos]{todonotes}
\usepackage{enumitem}
\usepackage[allcolors=blue]{hyperref}
\newtheorem{theorem}{Theorem}
\newtheorem{proposition}[theorem]{Proposition}
\newtheorem{lemma}[theorem]{Lemma}

\newtheorem{remark}[theorem]{Remark}

\theoremstyle{definition}
\newcommand{\R}{\mathbb{R}}

\newcommand{\Sf}{\mathbb{S}}
\newcommand{\C}{\mathbb{C}}

\newcommand{\spa}{\mbox{span}}
\newcommand{\hess}{\mbox{Hess\,}}

\newcommand{\nab}{\tilde\nabla}

\newcommand{\End}{\mbox{End}}
\newcommand{\Hom}{\mbox{Hom}}
\newcommand{\trace}{\mbox{tr\,}}
\renewcommand{\o}{\omega}

\newcommand{\Les}{\mathbb{L}}
\def\<{{\langle}}
\def\>{{\rangle}}

\def\Sal{{\cal S}}

\def\T{{\cal T}}

\def\B{{\cal B}}
\def\n{\nabla}
\def\d{\partial}
\def\a{\alpha}

\def\e{\epsilon}
\def\be{\begin{equation} }
\def\ee{\end{equation} }

\def\pu{\partial_u }
\def\pv{\partial_v }
\def\pz{\partial_z }
\renewcommand{\gg}{\Gamma^1 }
\newcommand{\gh}{\Gamma^2 }
\def\proof{\noindent{\it Proof:  }}
\def\qed{\ifhmode\unskip\nobreak\fi\ifmmode\ifinner
\else\hskip5 pt \fi\fi\hbox{\hskip5 pt \vrule width4 pt
height6 pt  depth1.5 pt \hskip 1pt }}
\begin{document}

\title{Conformal infinitesimal variations of\\
Euclidean hypersurfaces}
\author{M. Dajczer, M. I. Jimenez and Th. Vlachos}
\date{}
\maketitle

\begin{abstract} In the realm of conformal geometry, we give a 
parametric classification of the hypersurfaces in Euclidean space 
that admit nontrivial conformal infinitesimal variations. A 
parametric classification of the Euclidean hypersurfaces that allow 
a nontrivial conformal variation was obtained by E. Cartan in 1917. 
In particular, we show that the class of hypersurfaces studied
here is much larger than the one characterized by Cartan.
\end{abstract}

Classifying Euclidean hypersurfaces $f\colon M^n\to\R^{n+1}$, 
$n\geq 3$, that admit smooth isometric or, more generally, 
conformal variations is a classical subject in differential 
geometry. It was considered, among others, by U. Sbrana \cite{Sb1} 
and E. Cartan \cite{Ca1}, \cite{Ca2}
during the first part of the $20^{th}$ century. The case 
of the variations that are only infinitesimally isometric was 
studied by Sbrana \cite{Sb2}.
Modern  presentations of their results, as well as a large amount 
of additional information, can be found in \cite{DJ1}, \cite{DJ2}, 
\cite{DJ3}, \cite{DT1}, \cite{DT2}, \cite{DV} and \cite{Ji}.

Cartan \cite{Ca2} in a rather long and very difficult paper gave a local 
parametric classification of the hypersurfaces $M^n$ in $\R^{n+1}$ 
of dimension $n\geq 5$ that admit nontrivial conformal variations. 
The case $n=4$ was subsequently treated by Cartan \cite{Ca3} but 
only for a special class of submanifolds, thus the full classification 
for this dimension remains an open problem. The hypersurfaces in 
Cartan's description are conformally surface-like, conformally 
flat, conformally ruled and certain two-parameter congruences of 
hyperspheres. That $f\colon M^n\to\R^{n+1}$ is \emph{conformally 
surface-like} means that it is conformally congruent 
to either a cylinder or a rotation hypersurface over a surface 
in $\R^3$ or a cylinder over a three-dimensional hypersurface 
of $\R^4$  that is a cone over a surface in the sphere 
$\Sf^3\subset\R^4$. The hypersurface is \emph{conformally ruled} 
if $M^n$ carries an integrable $(n-1)$-dimensional distribution 
such that the restriction of $f$ to any leaf is an umbilical 
submanifold of $\R^{n+1}$.  Cartan in the same paper also proved that
conformally flat hypersurfaces are characterized by possessing 
at any point a principal curvature of multiplicity at least 
$n-1$. Hence, the really interesting class in Cartan's 
classification, and target of much of his work, is the one that
consists of certain two-parameter congruences of hyperspheres and are
parametrized in such way. These are nothing else that hypersurfaces 
that have a principal curvature of constant multiplicity $n-2$.

A \emph{conformal variation} of a given Euclidean hypersurface 
$f\colon M^n\to\R^{n+1}$ is a smooth variation 
$F\colon I\times M^n\to\R^{n+1}$, with $0\in I\subset\R$ 
an open interval and $F(0,\cdot)=f$, such that $f_t=F(t,\cdot)$ 
is a conformal immersion with respect to the metric induced by 
$f$ for any $t\in I$. The variation is said to be \emph{trivial} 
if it is induced by a composition of $f$ with a 
smooth family of conformal transformations of the Euclidean 
ambient space; recall that the latter transformations
are characterized by Liouville's theorem. 

If $F$ is a conformal variation of $f\colon M^n\to\R^{n+1}$, 
then there is a positive function 
$\gamma\in C^\infty(I\times M)$ with $\gamma(0,x)=1$ such that 
$$
\gamma(t,x)\<f_{t*}X,f_{t*}Y\>=\<X,Y\>
$$
for any tangent vector fields $X,Y\in\mathfrak{X}(M)$. The 
derivative of this equation with respect to $t$ computed 
at $t=0$ yields that the variational vector field 
$\T=F_*\d/\d t|_{t=0}$ of $F$ has to satisfy the condition
$$
\<\nab_X\T,f_*Y\>+\<f_*X,\nab_Y\T\>=2\rho\<X,Y\>
$$
for any vector fields $X,Y\in\mathfrak{X}(M)$, where 
the function $\rho\in C^\infty(M)$ 
given by $\rho(x)=-(1/2)\d\gamma/\d t(0,x)$ is called 
the conformal factor of $\T$.  Here and elsewhere we 
use the same notation for the inner products in $M^n$ and 
$\R^{n+1}$. Moreover, we denote by $\n$ and $\nab$, 
respectively, the  Levi-Civita connections associated to 
the induced metric in $M^n$ and the flat metric of the 
ambient space. 
\vspace{1ex}

The purpose of this paper is to classify parametrically the
hypersurfaces in Euclidean space $f\colon M^n\to\R^{n+1}$, 
$n\geq 5$, that  allow conformal infinitesimal variations. 
These variations are just the infinitesimal analogue to 
the conformal variations discussed above. 
\vspace{1ex}

A \emph{conformal infinitesimal variation} of a given Euclidean 
hypersurface $f\colon M^n\to\R^{n+1}$ is a smooth variation 
$F\colon I\times M^n\to\R^{n+1}$ of $f$ such that there is a
function $\gamma\in C^\infty(I\times M)$ satisfying 
$\gamma(0,x)=1$ and 
\be\label{varcond}
\frac{\d}{\d t}|_{t=0}\left(\gamma(t,x)\<f_{t*}X,f_{t*}Y\>\right)=0
\ee
for any vector fields $X,Y\in\mathfrak{X}(M)$.
\vspace{1ex} 

It has been shown in \cite{DJ2} that the study of conformal 
infinitesimal variations belongs to the realm of conformal 
geometry. In fact, by composing such a  variation with a
conformal transformation of the ambient Euclidean space
we obtain a new variation in the same class.
\vspace{1ex} 

It follows from \eqref{varcond} that the variational vector
field of a conformal infinitesimal variation satisfies the
same condition than the one of a conformal variation. This 
leads to the following definition.
\vspace{1ex}

A \emph{conformal infinitesimal bending} with \emph{conformal 
factor} $\rho\in C^\infty(M)$ of an isometric immersion 
$f\colon M^n\to\R^{n+1}$ of a Riemannian manifold $M^n$ 
into Euclidean space is a smooth section 
$\T\in\Gamma(f^*T\R^{n+1})$ that satisfies the condition
\be\label{cib}
\<\nab_X\T,f_*Y\>+\<f_*X,\nab_Y\T\>=2\rho\<X,Y\>
\ee
for any  $X,Y\in\mathfrak{X}(M)$.
\vspace{1ex}

A conformal infinitesimal bending is called a \emph{trivial 
conformal infinitesimal bending} 
if it is locally the restriction of a conformal Killing vector 
field of the Euclidean ambient space to the hypersurface. 
It is well-known that any conformal Killing field on an open 
connected subset of $\R^n$, $n\geq 3$, has the form
$$
\chi(x)=(\<x,v\>+\lambda)x-(1/2)\|x\|^2v+Cx+w
$$
where  $\lambda\in\R$, $v,w\in\R^n$,  $C\in\End(\R^n)$
is skew-symmetric and the conformal factor is 
$\rho=\<x,v\>+\lambda$; cf.\ \cite{Sc} for details.
\vspace{1ex}

It is already known from classical differential geometry  
that the appropriate approach to study infinitesimal variations 
is to deal with the variational vector field. In our case,
this statement is based in the following observation that establishes 
a kind of equivalence for the existence of nontrivial conformal 
infinitesimal bendings and variations. 
\vspace{1ex}

\noindent\emph{Correspondence:} Given an isometric immersion 
$f\colon M^n\to\R^{n+1}$, associated to any conformal 
infinitesimal variation of $f$ there is a conformal infinitesimal 
bending. Conversely, if $\T$ is a conformal infinitesimal bending 
of $f$, then the smooth variation $F\colon\R\times M^n\to\R^{n+1}$ 
defined by 
\be\label{unique}
F(t,x)=f(x)+t\T(x)
\ee
is a conformal infinitesimal variation with variational 
vector field $\T$ since \eqref{varcond} is satisfied for 
$\gamma(t,x)=e^{-2t\rho(x)}$. By no means
\eqref{unique} is unique with this property although it may 
be seen as the simplest one. For instance, new variations with 
variational vector field $\T$ are obtained by adding to 
\eqref{unique} terms of the type $t^k\delta$ with $k>1$, where 
$\delta\in\Gamma(f^*T\R^{n+1})$ and, if necessary, for restricted
values of $t$.
\vspace{1ex}

In view of the above, we say that a conformal infinitesimal 
variation is a \emph{trivial conformal infinitesimal 
variation} if the associated conformal infinitesimal bending 
is trivial.
\vspace{1ex}

An \emph{infinitesimal bending} of $f\colon M^n\to\R^{n+1}$ 
is a conformal infinitesimal bending of $f$ whose conformal 
factor is $\rho=0$. It is called \emph{trivial} if it is locally 
the restriction to the hypersurface of a Killing vector field of 
$\R^{n+1}$. Let $\T_1$ be a conformal infinitesimal bending of 
$f$ with conformal factor $\rho$ and let $\T_0$ be an infinitesimal 
bending of $f$. Then $\T_2=\T_1+\T_0$ satisfies \eqref{cib}, 
and thus it is also a conformal infinitesimal bending of $f$ with 
conformal factor $\rho$. In this paper, we always identify 
two conformal infinitesimal bendings of $f$ with the 
same conformal factor if they differ by a trivial infinitesimal 
bending. We also identify a conformal infinitesimal bending 
$\T$ with any of its constant multiples $c\,\T$, $0\neq c\in\R$.
\vspace{1ex}

Before we state our results, we recall from \cite{DJ2} that in 
order to admit a nontrivial conformal infinitesimal variation 
the hypersurface must possess a principal curvature of multiplicity 
at least $n-2$ at any point.
\vspace{1ex}

We first deal with the interesting class that includes, 
but is \emph{much larger} than, the interesting class of 
two-parameter congruences of hyperspheres in Cartan's 
classification of hypersurfaces that admit conformal 
variations; see Remark \ref{cartan}.  

\begin{theorem}\label{main} Let $f\colon M^n\to\R^{n+1}$, 
$n\geq 5$, admit a nontrivial conformal infinitesimal 
variation.  Assume that $f$ is neither conformally surface-like
nor conformally flat nor conformally ruled on any open 
subset of $M^n$. Then, on each connected component of 
an open dense subset of $M^n$, the hypersurface can be 
parametrized in terms of the conformal Gauss parametrization 
by a special hyperbolic or a special elliptic pair.

Conversely, any hypersurface $f\colon M^n\to\R^{n+1}$, 
$n\geq 5$, given in terms of the conformal Gauss parametrization 
by a special hyperbolic or special elliptic pair 
admits a nontrivial conformal infinitesimal variation.
Moreover, the conformal infinitesimal bendings  associated to 
any pair of nontrivial conformal infinitesimal variations of $f$ 
differ by a trivial conformal infinitesimal bending.
\end{theorem}

Special hyperbolic and elliptic pairs are the object  
of next section. As for the conformal Gauss parametrization, 
it goes as follows:  Let $f\colon\,M^n\to\R^{n+1}$, $n\geq 4$, be a 
hypersurface with Gauss map $N\colon M^n\to\Sf^n(1)\subset\R^{n+1}$ 
that possesses a principal curvature $\lambda>0$ of multiplicity $n-2$.
It is well-known that the corresponding eigenspaces form an 
integrable distribution and that $\lambda$ is constant along 
its leaves. Then, the focal map $f+rN\colon\,M^n\to\R^{n+1}$, 
where $r=1/\lambda$, induces an isometric immersion 
$h\colon\,L^2\to\R^{n+1}$. Here $L^2$ is the quotient space 
of leaves and $r\in C^{\infty}(L^2)$ satisfies $\|\nabla^h r\|<1$. 

Then $f$ can be locally parametrized along the unit normal 
bundle $N_1 L$ of $h$ by 
\be\label{Gausspar}
X(\xi)=h-r\left(h_*\nabla^h r+\sqrt{1-\|\nabla^h r\|^2}\,\xi\right).
\ee
Conversely, given a surface $h\colon\,L^2\to\R^{n+1}$ and 
$r\in C^{\infty}(L^2)$ positive whose gradient  satisfies 
$\|\nabla^h r\|<1$, then on the open subset of regular points, 
the parametrized hypersurface determined as in \eqref{Gausspar} 
by the pair $(h,r)$ has, with respect to the Gauss map
$N=h_*\nabla^h r+\sqrt{1-\|\nabla^h r\|^2}\,\xi$, the principal 
curvature $\lambda=1/r$ of multiplicity $n-2$.
\vspace{2ex}

We conclude this section with the case of ruled hypersurfaces.

\begin{theorem}\label{main2}
Let $f\colon M^n\to\R^{n+1}$, $n\geq 5$, be a conformally ruled 
hypersurface that is neither conformally surface-like nor 
conformally flat  on any open subset of $M^n$. Then, on each 
connected component of an open dense subset of $M^n$, $f$ admits  
a family of conformal infinitesimal bendings that are in one-to-one 
correspondence with the set of smooth functions on an interval. 
Moreover, any such bending is the variational vector field of a 
conformal variation.
\end{theorem}

\section{Special hyperbolic and elliptic pairs}

The purpose of this section is to introduce the notions 
of special hyperbolic and special elliptic pairs and 
show how they can be parametrically generated in terms
of a set of solutions of a second order hyperbolic or
elliptic PDE.
\vspace{2ex}

Let $g\colon L^2\to\Sf_1^{n+2}$ be a surface in the unit
Lorentzian sphere (de Sitter space) considered as a 
hypersurface of the Lorentzian space $\Les^{n+3}$, that is,
$$
\Sf^{n+2}_1=\{x\in\Les^{n+3}\colon\<x,x\>=1\}.
$$
We fix a pseudo-orthonormal basis $e_1,\ldots,e_{n+3}$ of 
$\Les^{n+3}$, that is,
$$
\|e_1\|=0=\|e_{n+3}\|, \<e_1,e_{n+3}\>=-1/2
\mbox{ and } \<e_i,e_j\>=\delta_{ij}\;\; \mbox{if}\;\; i\neq 1,n+3,
$$
and set
$g=(g_1,g_2,\ldots,g_{n+3})\colon L^2\to\Sf_1^{n+2}\subset\Les^{n+3}$
in terms of this basis. We assume that $g_1\neq 0$ everywhere, 
and let the map $h\colon L^2\to\R^{n+1}$ and the function
$r\in C^{\infty}(L)$ be given by 
\be\label{dere}
h=r(g_2,\ldots,g_{n+2})\;\;\mbox{and}\;\;r=1/g_1.
\ee
Notice that $g$ can be recovered from the pair $(h,r)$ by taking
\be\label{dere2}
g=r^{-1}(1,h,\|h\|^2-r^2).
\ee

\begin{proposition}\label{iff} We have that $L^2$ is a 
Riemannian manifold with the metric induced by $g$ if 
and only if $h$ is an immersion and the gradient of $r$ 
in the metric induced by $h$ satisfies $\|\nabla^h r\|<1$.
\end{proposition}

\proof See Lemma $12$ of \cite{DT1}.\vspace{2ex}\qed
 
A Riemannian surface $g\colon L^2\to\Sf_1^m$, 
$m\geq 4$, is said to be a \emph{hyperbolic surface} (respectively, 
\emph{elliptic surface}) if there is a tensor $I\neq J\in\End(TL)$ 
that verifies $J^2= I$ (respectively, $J^2=-I$) such that
the second fundamental form 
$\alpha^g\colon TL\times TL\to N_gL$
of $g$ satisfies
$$
\alpha^g(JX,Y)=\alpha^g(X,JY)
$$
for any $X,Y\in\mathfrak{X}(L)$. It is easily seen that
$J$ is unique up to sign. 

A local system of coordinates $(u,v)$ on $L^2$  is said 
to be \emph{real conjugate} for the surface $g\colon L^2\to\Sf_1^m$ 
if the condition
$$
\alpha^g(\d_u,\d_v)=0
$$
holds for the coordinate vector fields $\pu=\d/\d u$ and 
$\pv=\d/\d v$. The coordinate system is said to be 
\emph{complex conjugate} for $g$ if 
$$
\alpha^g(\d_z,\d_{\bar z})=0
$$
where $z=u+iv$ and $\pz=(1/2)(\pu-i\pv)$, that is, if 
$$
\alpha^g(\pu,\pu)+\alpha^g(\pv,\pv)=0.
$$

For a system of real conjugate coordinates let $\gg$, 
$\gh$ be the Christoffel symbols defined by 
\be\label{christsymbols}
\n_{\d_u}\d_v=\gg\d_u+\gh\d_v.
\ee
For a system of complex conjugate coordinates let 
$\Gamma$ be defined by 
\be\label{eq:gzbarz}
\n_{\d_z}\d_{\bar z}=\Gamma\d_z+\bar \Gamma\d_{\bar z},
\ee
where $\n$ also denotes the $\C$-bilinear extension of $\n$.
\vspace{2ex}

An elementary argument gives the following result.

\begin{proposition}\label{coordinates}  
If $g\colon L^2\to \Sf_1^m$ is a hyperbolic (respectively, 
elliptic) surface, then there exist local real conjugate 
(respectively, complex conjugate) coordinates 
on $L^2$ for $g$. Conversely, if there exist real conjugate 
(respectively, complex conjugate) coordinates on $L^2$, then 
$g\colon L^2\to\Sf_1^m$ is hyperbolic (respectively, 
elliptic). 
\end{proposition}

\proof See Proposition $11.10$ of \cite{DT1}.
\vspace{2ex}\qed

We call a hyperbolic surface  $g\colon L^2\to\Sf_1^m$ 
endowed with a system of real conjugate coordinates as in 
Proposition \ref{coordinates} a \emph{special hyperbolic 
surface} if the Christoffel symbols $\gg,\gh$ given by 
\eqref{christsymbols}  satisfy the condition
\be\label{integcon}
\gg_u=\gh_v.
\ee

\begin{proposition}\label{hypenv} Let $g\colon L^2\to\Sf_1^m$ 
be a simply connected special hyperbolic surface and let 
$\mu\in C^{\infty}(L)$ be the unique (up to a constant factor) 
positive solution of
\be\label{eq:siti2}
d\mu+2\mu\o=0
\ee
where $\o=\gh du+\gg dv$.
Then  $\varphi\in C^{\infty}(U)$ is a solution of the equation
\be\label{eq:eqvar}
\varphi_{uv}-\Gamma^1\varphi_u-\Gamma^2\varphi_v+F\varphi=0
\;\;\mbox{where}\;\;F=\<\d_u,\d_v\>
\ee
if and only if $\psi=\sqrt{\mu}\varphi$ satisfies 
\be\label{eq:eqpsi}
\psi_{uv}+M\psi=0
\ee
where
\be\label{Mm}
M=F-\frac{\mu_{uv}}{2\mu}+\frac{\mu_u\mu_v}{4\mu^2}\cdot
\ee
In particular, the map $k=\sqrt{\mu}\,h\colon L^2\to\Les^{m+1}$, 
where $h$ is the composition $h=i\circ g$ of $g$ with the 
inclusion $i\colon\Sf_1^m\to\Les^{m+1}$, satisfies
\be\label{form}
k_{uv}+Mk=0.
\ee

Conversely, for a system of coordinates $(u,v)$ on an open subset 
$U\subset\R^2$ let $\{k_1,\ldots,k_{m+1}\}$ be a set of 
solutions of \eqref{eq:eqpsi} for $M\in C^\infty(U)$. 
Assume that the map $k=(k_1,\ldots,k_{m+1})\colon U\to\Les^{m+1}$ 
satisfies $\mu=\|k\|^2>0$ and that the map
$h=(1/\sqrt{\mu})\,k\colon U\to\Les^{m+1}$ is an immersed 
surface with induced Riemannian metric. Then  
$g\colon U\to\Sf_1^m$ defined by  $h=i\circ g$ is a special 
hyperbolic surface.
\end{proposition}

\proof Notice that \eqref{integcon} is the integrability
condition of \eqref{eq:siti2}. Since $\mu\in C^\infty(L)$ is 
a solution of (\ref{eq:siti2}), it satisfies
$$
\gg=-\frac{\mu_v}{2\mu}\;\;\mbox{and}\;\;
\gh=-\frac{\mu_u}{2\mu}\cdot
$$
Hence (\ref{eq:eqvar}) becomes
\be\label{eq:guv2}
\varphi_{uv}+\frac{\mu_v}{2\mu}\varphi_u
+\frac{\mu_u}{2\mu}\varphi_v+F\varphi=0.
\ee
It follows easily that (\ref{eq:guv2}) takes the form 
(\ref{eq:eqpsi}) for $\psi=\sqrt{\mu}\,\varphi$, where 
$M$ is given by (\ref{Mm}).

We prove the converse. It is easily seen that 
$h=(1/\sqrt{\mu})\,k\colon U\to\Les^{m+1}$
satisfies
\be\label{eqh}
h_{uv}+\frac{\mu_v}{2\mu}h_u+\frac{\mu_u}{2\mu}h_v+Fh=0
\ee
where $F=M+\frac{\mu_{uv}}{2\mu}-\frac{\mu_u\mu_v}{4\mu^2}$.
If $h$ is an immersed Riemannian surface and 
$g\colon U\to\Sf_1^m$ is the surface defined by  $h=i\circ g$, 
then \eqref{eqh} implies that  $(u,v)$ are real conjugate 
coordinates for $g$ and that the Christoffel symbols of 
the metric induced by $g$ are 
$$
\Gamma^1=-\frac{\mu_v}{2\mu}
\;\;\mbox{and}\;\;\Gamma^2=-\frac{\mu_u}{2\mu}\cdot
$$
It follows that \eqref{integcon} is satisfied and that $\mu$ 
is a positive solution of \eqref{eq:siti2}.\vspace{2ex}\qed

We call an elliptic surface  $g\colon L^2\to\Sf_1^m$ endowed 
with a system of complex conjugate coordinates as in 
Proposition \ref{coordinates} a \emph{special elliptic surface} 
if the Christoffel symbol $\Gamma$ given by 
\eqref{eq:gzbarz}  satisfies the condition 
\be\label{integcon2}
\Gamma_z=\bar{\Gamma}_{\bar z}, 
\ee
that is, $\Gamma_z$ is real-valued. 

\begin{proposition}\label{ellenv} Let $g\colon L^2\to\Sf_1^m$ be a 
simply connected special elliptic surface and let 
$\mu\in C^{\infty}(L)$ be the unique (up to a constant factor) 
real-valued positive solution of
\be\label{eq:comp2}
\mu_{\bar{z}}+2\mu\Gamma=0.
\ee
Then $\varphi\in C^{\infty}(L)$ is a solution of 
\be\label{eq:eqvarc}
\varphi_{z\bar z}-\Gamma\varphi_z-\bar\Gamma\varphi_{\bar z}
+F\varphi=0\;\;\mbox{where}\;\;
F=\<\d_z,\d_{\bar z}\>=(1/4)(\|\pu\|^2+\|\pv\|^2)
\ee
if and only if $\psi=\sqrt{\mu}\varphi$ satisfies 
\be\label{comp}
\psi_{z\bar z}+M\psi=0,
\ee
where  
\be\label{Mm2}
M=F-\frac{\mu_{z\bar z}}{2\mu}+\frac{\mu_z\mu_{\bar z}}{4\mu^2}\cdot
\ee
In particular, the map
$k=\sqrt{\mu}\,h\colon L^2\to\Les^{m+1}$, where $h$ is the 
composition $h=i\circ g$ of $g$ with the inclusion 
$i\colon\Sf_1^m\to\Les^{m+1}$, satisfies
\be\label{form2}
k_{z\bar z}+Mk=0.
\ee

Conversely, for a system of coordinates $(u,v)$ on an open subset 
$U\subset\R^2$ let $\{k_1,\ldots,k_{m+1}\}$ be a set of 
solutions of  \eqref{comp} for $M\in C^\infty(U)$. 
Assume that the map $k=(k_1,\ldots,k_{m+1})\colon U\to\Les^{m+1}$ 
satisfies $\mu=\|k\|^2>0$ and that the map
$h=(1/\sqrt{\mu})\,k\colon U\to\Les^{m+1}$ is an immersed 
surface with induced Riemannian metric. Then  
$g\colon U\to\Sf_1^m$ defined by  $h=i\circ g$ is a special 
elliptic surface.
\end{proposition}

\proof Notice that \eqref{integcon2} is the integrability
condition of equation \eqref{eq:comp2}.
Since $\mu\in C^{\infty}(L)$ is a real-valued solution 
of (\ref{eq:comp2}) then $\Gamma=-(1/2\mu)\mu_{\bar z}$.
Hence (\ref{eq:eqvarc}) becomes
\be\label{eq:eqvarcm}
\varphi_{z\bar z}+\frac{\mu_z}{2\mu}\varphi_{\bar z}
+\frac{\mu_{\bar z}}{2\mu}\varphi_z+F\varphi=0.
\ee
It follows easily that (\ref{eq:eqvarcm}) has the 
form (\ref{form2}) for $k=\sqrt{\mu}\,\varphi$ where 
$M$ is given by~(\ref{Mm2}).

We prove the converse.  It is easily seen that 
$h=(1/\sqrt{\mu})\,k\colon U\to\Les^{m+1}$
satisfies
\be\label{eqh2}
h_{z\bar z}+\frac{\mu_z}{2\mu}h_{\bar z}
+\frac{\mu_{\bar z}}{2\mu}h_z+Fh=0
\ee
where $F=M+\frac{\mu_{z\bar z}}{2\mu}
-\frac{\mu_z\mu_{\bar z}}{4\mu^2}$.
If $h$ is an immersed Riemannian surface and 
$g\colon U\to\Sf_\e^m$ is the surface defined by $h=i\circ g$, 
then \eqref{eqh2} implies that $(u,v)$ are complex 
conjugate coordinates for $g$ and that the complex Christoffel 
symbol of the metric induced by $g$ is 
$\Gamma=-(1/2\mu)\mu_{\bar z}$.
It follows that \eqref{integcon2} is satisfied and that  
$\mu$ is a  positive solution of \eqref{eq:comp2}.
\vspace{2ex}\qed

We call the pair $(h,r)$ formed by a surface 
$h\colon L^2\to\R^m$ and a function $r\in C^{\infty}(L)$
a \emph{special hyperbolic pair} (respectively, 
\emph{special elliptic pair}) if there exists a special hyperbolic 
surface (respectively, special elliptic surface) 
$g\colon L^2\to\Sf_1^{m+1}$ such that $(h,r)$ are given by 
\eqref{dere}.\vspace{2ex}

We conclude this section with the following result that 
will be used for the proof of Theorem \ref{main}.

\begin{proposition}\label{equivspecial} For a simply 
connected surface $g\colon L^2\to\Sf_1^m\subset\Les^{m+1}$ 
the following assertions are equivalent:
\begin{itemize}
\item[(i)] The surface $g$ is special hyperbolic 
(respectively, special elliptic).
\item[(ii)] The surface is hyperbolic (respectively, elliptic) 
with respect to a tensor $J$ on $L^2$ satisfying $J^2=I$ 
(respectively, $J^2=-I$) and there exists 
$\mu\in C^\infty(L)$ nowhere vanishing such that $D=\mu J$ 
is a Codazzi tensor on $L^2$.
\end{itemize}
\end{proposition}

\proof  Let $g$ be a hyperbolic surface as in part $(ii)$ 
and let $(u,v)$ be local real conjugate  coordinates on $L^2$ 
given by Proposition \ref{coordinates}.  Then the equation 
\be\label{eq:codbard}
\left(\n_{\pu}D\right)\pv
-\left(\n_{\pv}D\right)\pu=0
\ee
is easily seen to be equivalent to (\ref{eq:siti2}). 

Conversely, if $g$ is special hyperbolic with real conjugate 
coordinates $(u,v)$, $J\in\Gamma(\End(T))$ is 
given by $J\pu=\pu$ and $J\pv=-\pv$, and $\mu\in C^{\infty}(L)$ 
satisfies (\ref{eq:siti2}), then $D=\mu J$ satisfies 
(\ref{eq:codbard}) in view of (\ref{eq:siti2}), and hence 
is a Codazzi tensor on $L^2$.  The proof for the elliptic 
case is similar.\qed

\section{The proofs}

Let $f\colon M^n\to\R^{n+1}$ be a hypersurface with second 
fundamental form $A$ with respect to the Gauss map 
$N\in\Gamma(N_fM)$.
Associated to a conformal infinitesimal bending $\T$ with 
conformal factor $\rho$ there is a symmetric tensor 
${\cal B}\in\Gamma(\End(TM))$ defined as follows:
Let $L\in\Gamma(\Hom(TM,f^*T\R^{n+1}))$ be the tensor
defined by
$$
LX=\nab_X \T-\rho f_*X=\T_*X-\rho f_*X
$$
for any $X\in\mathfrak{X}(M)$. Then let
$B\colon TM\times TM\to f^*T\R^{n+1}$ be given by
$$
B(X,Y)=(\nab_XL)Y=\nab_XLY-L\n_XY
$$
for any $X,Y\in\mathfrak{X}(M)$. 
We define ${\cal B}\in\Gamma(\End(TM))$ by
$$
\<\B X,Y\>=\<B(X,Y),N\>
$$
for any $X,Y\in\mathfrak{X}(M)$. Notice that flatness
of the ambient space and 
$$
B(X,Y)_{N_fM}
=(\nab_X\nab_Y\T-\nab_{\nabla_XY}\T)_{N_fM}-\rho\<AX,Y\>N
$$
give that $\B$ is symmetric.

By Proposition $5$ in \cite{DJ2} the tensor $\B$ 
satisfies  the following fundamental system of equations:
\be\label{Gauss}
{\cal B}X\wedge AY-{\cal B}Y\wedge AX
+X\wedge HY-Y\wedge HX=0
\ee
and
\be\label{Codazzi}
(\nabla_X {\cal B})Y-(\nabla_Y {\cal B})X
+(X\wedge Y)A\nabla\rho=0
\ee
where $H\in\End(TM)$ is  defined by $HY=\nabla_Y\nabla\rho$
and $X,Y\in\mathfrak{X}(M)$. The equations are called  
fundamental because they are the integrability condition 
for the system of equations that determines a conformal 
infinitesimal bending; see Corollary~$2$ in \cite{DJ2}.
\vspace{1ex}

Trivial 
conformal infinitesimal bendings can be characterized in terms 
of their associated tensors as follows.

\begin{proposition}\label{trivial}
A conformal infinitesimal bending $\T$ of 
$f\colon M^n\to \R^{n+1}$, $n\geq3$, is trivial if and only 
if its associated tensor $\B$ has the form
$\B=\varphi I$ for $\varphi\in C^\infty(M)$.
\end{proposition}
\proof
See Corollary 9 in \cite{DJ2}.
\qed
\vspace{2ex}

Let $\T$ be a conformal infinitesimal bending of 
$f\colon M^n\to \R^{n+1}$ with conformal factor $\rho$. 
At any point of $M^n$ we obtain from \eqref{Gauss} that 
the associated bilinear form $\theta\colon TM\times TM\to\R^4$ 
defined by
{\em \be\label{theta}
\!\theta(X,Y)\!=\!(\<(A+\B)X,Y\>,\<(I+H)X,Y\>,\<(A-\B)X,Y\>,
\<(I-H)X,Y\>)
\ee}
is flat  with respect to the inner product $\<\!\<\,,\,\>\!\>$
of signature $(1,1,-1,-1)$. That $\theta$ is \emph{flat} means 
that 
$$
\<\!\<\theta(X,Y),\theta(Z,W)\>\!\>
-\<\!\<\theta(X,W),\theta(Z,Y)\>\!\>=0
$$
for all $X,Y,Z,W\in \mathfrak{X}(M)$.

\begin{proposition}\label{nulltriv}
Let $f\colon M^n\to \R^{n+1}$, $n\geq 3$, be an 
isometric immersion with no umbilical points. If 
$\T$ is a conformal infinitesimal bending of $f$ 
such that the associated flat bilinear form $\theta$ 
given by \eqref{theta} is null at any point of $M^n$
then $\T$ is trivial. 
\end{proposition}

\proof That $\theta$ is \emph{null} means that
$$
\<\!\<\theta(X,Y),\theta(Z,W)\>\!\>=0
$$
for all $X,Y,Z,W\in \mathfrak{X}(M)$. Equivalently,
$$
\<AX,Y\>\B+\<\B X,Y\>A+\<X,Y\>H+\<HX,Y\>I=0
$$
for any $X,Y\in\mathfrak{X}(M)$. 
Fix a point $x\in M^n$. By the above $A(x),\B(x)$ 
and $H(x)$ commute, that is, there exists  an 
orthonormal basis $\{X_i\}_{1\leq i\leq n}$ of $T_xM$ 
that diagonalizes them simultaneously. 
If $\lambda_i$, $b_i$ and $h_i$ are the respective 
eigenvalues of $A(x),\B(x)$ and $H(x)$ corresponding 
to $X_i$, $1\leq i\leq n$, then 
$$
\lambda_i\B+b_iA+H+h_iI=0.
$$
Let $i,j$ be such that $\lambda_i\neq\lambda_j$. Then
$$
(\lambda_i-\lambda_j)\B+(b_i-b_j)A+(h_i-h_j)I=0.
$$
Hence 
$$
(\lambda_i-\lambda_j)b_i+(b_i-b_j)\lambda_i+(h_i-h_j)
=(\lambda_i-\lambda_j)b_j+(b_i-b_j)\lambda_j+(h_i-h_j)=0,
$$
and therefore
$$
(\lambda_i-\lambda_j)(b_i-b_j)=0
$$
showing that $b_i=b_j$.
If $i\neq j$ are such that $\lambda_i=\lambda_j$ then 
$$
(b_i-b_j)A+(h_i-h_j)I=0.
$$
But since $f$ has no umbilical points, we necessarily have 
$b_i=b_j$ and hence $\B=bI$ at any $x\in M^n$. We conclude 
from Proposition \ref{trivial} that $\T$ is trivial.\qed

\begin{lemma}\label{eigen}
Let $f\colon M^n\to\R^{n+1}$ be an isometric immersion free
of umbilical points and let $\T$  be a nontrivial conformal
infinitesimal bending of $f$. Then $A$, $\B$ and $H$ share 
a common eigenbundle $\Delta$ with $\dim\Delta\geq n-2$
on connected components of an open and dense subset of $M^n$.
\end{lemma}

\proof By  Proposition \ref{nulltriv}
the bilinear form $\theta$ is not null. Theorem $3$
in \cite{DF}  or Lemma $4.22$ in \cite{DT2} yield an
orthogonal decomposition
$\R^4=\R^{2,2}=\R^{\ell,\ell}\oplus\R^{2-\ell,2-\ell},
\;1\leq\ell\leq 2$, such that  the $\R^{\ell,\ell}$-component
$\theta_1$ of $\theta$ is nonzero but is null since 
$\Sal(\theta_1)=\Sal(\theta)\cap\Sal(\theta)^\perp$, and the
$\R^{2-\ell,2-\ell}$-component $\theta_2$ is flat and satisfies
$\dim\mathcal{N}(\theta_2)\geq n-4+2\ell$. Moreover, since
$\theta$ is not null then  $\ell=1$.

We denote $\Delta=\mathcal{N}(\theta_2)$ and restrict ourselves
to connected components of an open and dense subset where
$\dim\Delta\geq n-2$ is constant. Since we have that
$\theta(T,X)=\theta_1(T,X)$ for any $T\in\Gamma(\Delta)$
and $X\in\mathfrak{X}(M)$, then
$$
\<\!\<\theta(T,X),\theta(Y,Z)\>\!\>=0
$$
for any $T\in\Gamma(\Delta)$ and $X,Y,Z\in\mathfrak{X}(M)$.
Equivalently,
\be\label{ABHI}
\<AT,X\>\B+\<\B T,X\>A+\<T,X\>H+\<HT,X\>I=0
\ee
for any $T\in\Gamma(\Delta)$ and $X\in\mathfrak{X}(M)$.
Taking $X$ orthogonal to $T$ we see that
\be\label{ABI}
\<AT,X\>\B+\<\B T,X\>A+\<HT,X\>I=0.
\ee
Fix $x\in M^n$ and assume that there exists $T\in\Delta(x)$
and $X\in T_xM$ such that $\<X,T\>=0$ and $\<\B T,X\>\neq 0$.
From \eqref{ABI} and since $f$ is free of umbilic points we
have that $A$ commutes with $\B$. Hence also does $H$.
Let $\{X_i\}_{1\leq i\leq n}$ be an orthonormal basis of $T_xM$
of common eigenvectors of $A$, $\B$ and $H$ with corresponding
eigenvalues $\lambda_i$, $b_i$ and $h_i$. Since $\<\B T,X\>\neq 0$
with $\<X,T\>=0$, then $T$ is not an eigenvector. Hence, there are
two eigenvalues $b_1\neq b_2$ such that 
$\<T,X_1\>\neq 0\neq\<T,X_2\>$.  
Thus, we have from \eqref{ABHI} that
$$
\lambda_1\B+b_1A+H+h_1I=0\;\;\mbox{and}
\;\;\lambda_2\B+b_2A+H+h_2I=0.
$$
Hence
\be\label{ABI2}
(\lambda_1-\lambda_2)\B+(b_1-b_2)A+(h_1-h_2)I=0,
\ee
from where we obtain
$$
(\lambda_1-\lambda_2)b_j+(b_1-b_2)\lambda_j
+h_1-h_2=0,\;\;1\leq j\leq n.
$$
Taking the difference between the cases $j=1$ and $j=2$ we have
$$
(\lambda_1-\lambda_2)(b_1-b_2)=0,
$$
and therefore $\lambda_1=\lambda_2$. It follows from \eqref{ABI2}
that $A$ is a multiple of the identity which is a contradiction.

Therefore $\<\B T,X\>=0$ for any $T\in\Delta(x)$ and $X\in T_xM$
with $\<X,T\>=0$. This implies that $\Delta$ is an eigenspace of
$\B$. If  $\<AT,X\>\neq 0$, for some $T\in\Delta(x)$ and
$X\in T_xM$ with $\<T,X\>=0$, then we have from \eqref{ABI} that
$B$ is a multiple of the identity and this is contradiction.
Hence $\Delta$ is also an eigenspace of $A$, and consequently of
$H$.
\vspace{2ex}\qed

Let  $f\colon M^n\to\R^{n+1}$ be a hypersurface that carries a  
principal curvature of constant multiplicity $n-2$ with 
corresponding eigenbundle $\Delta$. Recall that the 
\emph{splitting tensor} 
$C\colon\Gamma(\Delta)\to\Gamma(\End(\Delta^\perp))$ of $\Delta$
is defined by 
$$
C_TX=-\nabla_X^hT=-(\nabla_XT)_{\Delta^\perp}
$$ 
for any $T\in\Gamma(\Delta)$ and $X\in\Gamma(\Delta^\perp)$.
If $f$ is not conformally surface-like on any open subset of 
$M^n$ we say that $f$ is \emph{hyperbolic} (respectively,
\emph{parabolic} or \emph{elliptic}) if there exists 
$J\in\Gamma(\End(\Delta^\perp))$ satisfying the 
following conditions:
\begin{itemize}
\item[(i)] $J^2=I$ and $J\neq I$ (respectively, $J^2=0$, with $J\neq 0$, 
and $J^2=-I$),
\item[(ii)] $\nabla^h_T J=0$ for all $T\in\Gamma(\Delta)$,
\item[(iii)] $C_T\in\spa\{I,J\}$ for all $T\in\Gamma(\Delta)$.
\end{itemize}

Let $\mathbb{V}^{n+2}\subset\Les^{n+3}$ be the 
\emph{light cone} of $\Les^{n+3}$, that is,
$$
\mathbb{V}^{n+2}=\{v\in\Les^{n+3}\colon\<v,v\>=0,v\neq 0\}.
$$
Given $w\in \mathbb{V}^{n+2}$ we have that
$$
\mathbb{E}^{n+1}=\{v\in\mathbb{V}^{n+2}\colon \<v,w\>=1\}
$$
is a model of $\R^{n+1}$ in $\Les^{n+3}$. In fact, 
fix $v\in\mathbb{E}^{n+1}$ and a linear isometry 
$C\colon \R^{n+1}\to (\spa\{v,w\})^{\perp}\subset \Les^{n+3}$.
The  map $\Psi\colon\R^{n+1}\to\mathbb{V}^{n+2}\subset\Les^{n+3}$ 
given~by
\be\label{psi}
\Psi(x)=v+Cx-\frac{1}{2}\|x\|^2w
\ee
is an isometric embedding such that 
$\Psi(\R^{n+1})=\mathbb{E}^{n+1}$.

The normal bundle of $\Psi$ is $N_\Psi\R^{n+1}=\spa\{\Psi,w\}$ 
and its second fundamental form is given by
\be\label{sff}
\a^{\Psi}(U,V)=-\<U,V\>w
\ee
for any $U,V\in T\R^{n+1}$. For further details we refer to 
Section $9.1$ in \cite{DT2}.

Let $f\colon M^n\to\R^{n+1}$, $n\geq 5$, be an oriented 
hypersurface with a principal curvature 
$\lambda$ of constant multiplicity $n-2$. By composing 
with an appropriate inversion, if necessary, 
and given that $f$ is orientable,  we can always assume 
that $\lambda>0$ at any point of $M^n$. Let $A$ be the 
second fundamental form  associated to the Gauss map $N$ 
of $f$ and let $\Delta(x)\subset T_xM$ be the eigenspace 
corresponding to $\lambda(x)$ at $x\in M^n$. Fix an embedding  
$\Psi$ as in \eqref{psi} and let $S\colon M^n\to\Les^{n+3}$ 
be the map given by
\be\label{S}
S(x)=\lambda(x)\Psi(f(x))+\Psi_*N(x).
\ee
Then $S(x)\in\Sf_1^{n+2}\subset\Les^{n+3}$ and
\be\label{S*}
S_*X=X(\lambda)\Psi(f(x))-\Psi_*f_*(A-\lambda I)X
\ee
for any $X\in\mathfrak{X}(M)$. From \eqref{S*} it follows 
that $S$ is constant along the leaves of $\Delta$. 
Let $L^2$ be the quotient space of leaves of $\Delta$ and 
let $\pi\colon M^n\to L^2$ be the canonical projection. 
Thus $S$ induces an immersion 
$s\colon L^2\to\Sf_{1}^{n+2}\subset\Les^{n+3}$ such that 
$S=s\circ\pi$.
Moreover, the metric $\<\,,\,\>'$ on $L^2$ induced by $s$ satisfies
\be\label{metricsur}
\<\bar{X},\bar{Y}\>'=\<(A-\lambda I)X,(A-\lambda I)Y\>
\ee
where $X,Y\in\mathfrak{X}(M)$ are the horizontal lifts of 
$\bar{X},\bar{Y}\in\mathfrak{X}(L^2)$.

\begin{proposition}\label{fundD}
Let $f\colon M^n\to \R^{n+1}$, $n\geq 5$, be an oriented 
hypersurface and let $\T$ be a nontrivial conformal infinitesimal 
bending of $f$. Assume that the principal curvature $\lambda$ 
of $A$ determined by $\Delta$ from Lemma \ref{eigen} is 
positive and has constant multiplicity $\dim\Delta=n-2$. 
Then, on each connected component of an open dense subset of $M^n$
either $f$ is conformally surface-like or is  hyperbolic, 
parabolic or elliptic with respect to 
$J\in\Gamma(\End(\Delta^\perp))$ and there exists 
$\mu\in\C^{\infty}(M)$ nowhere vanishing and constant along 
the leaves of $\Delta$ such that  
$D=\mu J\in\Gamma(\End(\Delta^\perp))$ satisfies:
\begin{itemize}
\item[\hypertarget{i}{(i)}] $(A-\lambda I)D$ is symmetric,
\item[\hypertarget{ii}{(ii)}] $\nabla_T^hD=0$,
\item[\hypertarget{iii}{(iii)}]$
(\nabla_X(A-\lambda I)D)Y
-(\nabla_Y(A-\lambda I)D)X=X\wedge Y(D^t\nabla\lambda),
$
\item[\hypertarget{iv}{(iv)}]
$\<(\nabla_YD)X
-(\nabla_XD)Y,\nabla\lambda\>+\hess\lambda(DX,Y)
-\hess\lambda(X,DY)$\\
$=\lambda(\<AX,(A-\lambda I)DY\>
-\<(A-\lambda I)DX,AY\>),$

\item[\hypertarget{v}{(v)}] 
$(A-\lambda I)DX\wedge(A-\lambda I)Y
-(A-\lambda I)DY\wedge(A-\lambda I) X=0$
\end{itemize}
for any $T\in\Gamma(\Delta)$ and $X,Y\in\Gamma(\Delta^\perp)$.

Conversely, assume that $f$ as above is either 
hyperbolic, parabolic or elliptic with respect to 
$J\in\Gamma(\End(\Delta^\perp))$ and there is  
$0\neq D=\mu J\in\Gamma(\End(\Delta^\perp))$ that satisfies 
conditions $(i)$ to $(v)$.  If $M^n$ is simply connected  
there exists a nontrivial conformal infinitesimal 
bending $\T$ of $f$ determined by $D$, unique up to 
trivial conformal infinitesimal bendings.
\end{proposition}

\proof We have from Lemma \ref{eigen} that $\Delta$ is a common 
eigenbundle for $A$, $\B$ and $H$. Thus $\B|_\Delta=bI$ and 
$H|_\Delta=hI$ where $b,h\in C^{\infty}(M)$. We obtain from 
\eqref{ABHI} that  
$$
bA+\lambda\B+H+hI=0.
$$
In particular $\lambda b+h=0$,
and thus locally 
\be\label{suma}
bA+\lambda(\B-bI)+H=0.
\ee

From \eqref{Codazzi} we have
\be\label{I}
T(b-\lambda\rho)=T(b)-\lambda T(\rho)=0
\ee
for any $T\in\Gamma(\Delta)$. Then \eqref{Codazzi} is 
equivalent to 
\be\label{cod}
(\nabla_X (\B-bI))Y-(\nabla_Y (\B-bI))X
+(X\wedge Y)(A\nabla\rho-\nabla b)=0.
\ee
It follows from \eqref{I} and \eqref{cod} that
\be\label{sin}
(\nabla^h_T(\B-bI))X=(\B-bI)C_TX
\ee
for any $X\in\Gamma(\Delta^\perp)$ and $T\in\Gamma(\Delta)$. 

We regard $A-\lambda I$ and $\B-b I$ as tensors on $\Delta^\perp$. 
We obtain from \eqref{sin} and the Codazzi equation 
$\nabla_T^hA=(A-\lambda I)C_T$ that
$$
(\B-bI)C_T=C_T^t(\B-bI)\;\;\mbox{and}\;\; 
(A-\lambda I)C_T=C_T^t(A-\lambda I).
$$
We have that $D\in\Gamma(\End(\Delta^\perp))$ defined by
$$
D=(A-\lambda I)^{-1}(\B-bI)
$$
satisfies $D\neq 0$ since $\T$ is nontrivial.
Hence
\begin{align*}
(A-\lambda I)DC_T&=(\B-bI)C_T
=C_T^t(\B-bI)= C_T^t(A-\lambda I)D\\
&=(A-\lambda I)C_TD,
\end{align*}
and therefore 
\be\label{CTD}
[D,C_T]=0.
\ee

We also have 
$$
(A-\lambda I)C_T D= (\nabla_T^hA)D
$$
and 
\begin{align*}
(A-\lambda I)DC_T&=(\B-b I)C_T
=\nabla_T^h(\B-b I)
=\nabla_T^h((A-\lambda I)D)\\
&=\nabla_T^h(AD)-\lambda\nabla_T^hD.
\end{align*}
Thus
$$
(A-\lambda I)\nabla_T^hD=(A-\lambda I)[D,C_T],
$$
and hence 
\be\label{Dpar}
\nabla_T^hD=0
\ee
for any $T\in\Gamma(\Delta)$.
 
It follows from \eqref{CTD}, \eqref{Dpar} and 
Corollary $11.7$ of \cite{DT2} that $D$ 
is projectable with respect to $\pi\colon M^n\to L^2$, 
that is, $D$ is the horizontal lift of a tensor $\bar D$ on 
$L^2$. Hence 
$$
\pi_*DX=\bar{D}\pi_*X=\bar{D}\bar{X}\circ\pi
\;\;\;\mbox{if}\;\;\;\pi_*X=\bar{X}\circ\pi.
$$
We have that \eqref{Gauss} reads as
$$
\B X\wedge AY-\B Y\wedge AX+X\wedge HY-Y\wedge HX=0.
$$
Since $H=\lambda(bI-\B)-bA$ from \eqref{suma}, then
\be\label{Gauss2}
(\B-bI)X\wedge(A-\lambda I)Y-(\B-bI)Y\wedge(A-\lambda I)X=0
\ee
for any $X,Y\in\mathfrak{X}(M)$. From \eqref{Gauss2} and the 
definition of $D$ we have 
$$
\!\<((A-\lambda I)DX\wedge(A-\lambda I)Y
\!-\!(A-\lambda I)DY\wedge(A-\lambda I)X)
(A-\lambda I)Z,(A-\lambda I)W\>\!=\!0
$$
for any $X,Y,Z,W\in\Gamma(\Delta^\perp)$. This implies that
$$
\<(\bar{D}\bar{X}\wedge\bar{Y}
-\bar{D}\bar{Y}\wedge\bar{X})\bar{Z},\bar{W}\>'=0
$$
for any $\bar{X}, \bar{Y}, \bar{Z}, \bar{W}\in\mathfrak{X}(L)$. 
In other words, we have
$$
\bar{D}\bar{X}\wedge\bar{Y}-\bar{D}\bar{Y}\wedge\bar{X}=0
$$
with respect to the metric $\<\,,\,\>'$. Thus $\trace\bar{D}=0$. 

We have that $\bar{D}$ has either two smooth distinct 
real eigenvalues, a single real eigenvalue of multiplicity 
two or a pair of smooth complex conjugate eigenvalues. 
Thus there is $\bar{\mu}\in C^\infty(L)$ such that 
$\bar{D}=\bar{\mu}\bar{J}$, $\bar{J}\neq I$, where the tensor 
$\bar{J}\in\Gamma(\End(TL))$ satisfies 
$\bar{J}^2=\epsilon I$, for $\epsilon= 1,0$ or $-1$.
Hence $D=\mu J$ where $J$ is the lifting of $\bar{J}$
and $\bar{\mu}=\mu\circ\pi$. In particular $\trace D=0$. 

If $\spa\{C_T:T\in\Delta\}\subset\spa\{I\}$ we have from 
Corollary 9.33 in \cite{DT2} that $f$ is conformally 
surface-like. Hence, we assume  
$\spa\{C_T:T\in\Delta\}\not\subset\spa\{I\}$ and obtain 
from \eqref{CTD} that $C_T\in\spa\{I,J\}$ for any 
$T\in\Gamma(\Delta)$.

We have from \eqref{Dpar} that
$$
T(\mu)J+\mu\nabla^h_TJ=0
$$
for any $T\in\Gamma(\Delta)$. Therefore 
$$
\epsilon T(\mu)I+\mu J\nabla^h_T J=0\;\;\mbox{and}\;\;
\epsilon T(\mu)I+\mu (\nabla^h_T J)J=0.
$$
Since $J^2=\epsilon I$ we obtain that $T(\mu)=0$, and hence 
$\nabla_T^h J=0$.
Thus, the hypersurface $f$ is either hyperbolic, parabolic or
elliptic. 

We have from \eqref{suma} that
\be\label{one}
X(b)AY+X(\lambda)\B Y-X(\lambda b)Y+b(\nabla_X A)Y
+\lambda(\nabla_X\B)Y+(\nabla_X H)Y=0.
\ee 
On the other hand, the Gauss equation yields
\be\label{two}
(\nabla_X H)Y-(\nabla_Y H)X=R(X,Y)\nabla\rho
=\<AY,\nabla\rho\>AX-\<AX,\nabla\rho\>AY.
\ee
Then \eqref{Codazzi}, \eqref{one}, \eqref{two} and 
the Codazzi equation imply that
\begin{align*}
X(b)AY&+X(\lambda)\B Y-X(\lambda b)Y
-Y(b)AX-Y(\lambda)\B X+Y(\lambda b)X\\
&-\lambda(X\wedge Y)A\nabla\rho
+\<AY,\nabla\rho\>AX-\<AX,\nabla\rho\>AY=0
\end{align*}
for any $X,Y\in\mathfrak{X}(M)$. Then
\begin{align*}
\<X,\nabla b-A\nabla\rho\>&(A-\lambda I)Y
-\<Y,\nabla b-A\nabla\rho\>(A-\lambda I)X\\
&+\<X,\nabla\lambda\>(\B-bI)Y
-\<Y,\nabla\lambda\>(\B-bI)X=0
\end{align*}
for any $X,Y\in\mathfrak{X}(M)$.
For $X,Y\in\Gamma(\Delta^\perp)$ we obtain 
$$
\<X,\nabla b-A\nabla\rho\>Y-\<Y,\nabla b-A\nabla\rho\>X
+\<X,\nabla\lambda\>DY-\<Y,\nabla\lambda\>DX=0.
$$
Taking $X$ and $Y$ orthonormal, we obtain
$$
\<Y,\nabla b-A\nabla\rho\>-\<X,\nabla\lambda\>\<DY,X\>
+\<Y,\nabla\lambda\>\<DX,X\>=0
$$
and
$$
\<X,\nabla b-A\nabla\rho\>+\<X,\nabla\lambda\>\<DY,Y\>
-\<Y,\nabla\lambda\>\<DX,Y\>=0.
$$
Using that $\trace D=0$ this gives
\be\label{gradients}
D^t\nabla\lambda=\nabla b-A\nabla\rho
\ee
where $D^t$ denotes the transpose of $D$.

So far we have that $\hyperlink{i}{(i)}$ holds from 
the definition of $D$, $\hyperlink{ii}{(ii)}$ is 
\eqref{Dpar}, $\hyperlink{iii}{(iii)}$ follows from 
\eqref{cod} and \eqref{gradients}, and 
$\hyperlink{v}{(v)}$ is \eqref{Gauss2}. It remains 
to prove that $\hyperlink{iv}{(iv)}$ holds. 
To do this, fix a pseudo-orthonormal basis $e_1\dots,e_{n+3}$ 
of $\Les^{n+3}$ and set $v=e_1$ and $w=-2e_{n+3}$.
Let $\Psi\colon\R^{n+1}\to\mathbb{V}^{n+2}\subset\Les^{n+3}$ 
and $S\colon M^n\to\Les^{n+3}$ be given by \eqref{psi} 
and \eqref{S} respectively. We see next that the 
immersion $s\colon L^2\to\Sf_1^{n+2}\subset\Les^{n+3}$ 
induced by $S$ satisfies $s=g$, where $g$ is given by 
\eqref{dere2}, $h\colon L^2\to \R^{n+1}$ is induced 
by $f+rN$ and $r=\lambda^{-1}$. 
In fact, we have that $\Psi(y)=(1,y,\|y\|^2)$. Then
\begin{align*} 
S(x)&=\lambda(1,f(x),\|f(x)\|^2)+(0,N(x),2\<f(x),N(x)\>)\\
&=\lambda(1,f(x)+rN,\|f(x)\|^2+2r\<f(x),N(x)\>).
\end{align*}
Since $h\circ\pi=f+rN$, it follows that 
$$
s=r^{-1}(1,h,\|h\|^2-r^2)=g.
$$ 
Let $X,Y\in\Gamma(\Delta^\perp)$ be the horizontal lifts 
of $\bar{X},\bar{Y}\in\mathfrak{X}(L)$. We have
\begin{align*}
\nab'_XS_*DY&=\nab'_{\pi_*X}g_*\pi_*DY
=\nab'_{\bar{X}}g_*\bar{D}\bar{Y}\\
&=g_*\nabla'_{\bar{X}}\bar{D}\bar{Y}
+\a^g(\bar{X},\bar{D}\bar{Y})
-\<\bar{X},\bar{D}\bar{Y}\>'g\circ\pi
\end{align*}
where $\nab'$ and $\nabla'$ denote the connections in $\Les^{n+3}$
and $L^2$, respectively. We obtain from \eqref{S*} that
\begin{align*}
&\nab_X'\Psi_*f_*(A-\lambda I)DY
=X\<DY,\nabla\lambda\>\Psi\circ f
+\<DY,\nabla\lambda\>\Psi_*f_*X\\
&-g_*\nabla'_{\bar{X}}\bar{D}\bar{Y}
-\a^g(\bar{X},\bar{D}\bar{Y})
+\<(A-\lambda I)X,(A-\lambda I)DY\>
(\lambda\Psi\circ f+\Psi_*N).
\end{align*}
On the other hand, using \eqref{sff} and 
\eqref{S*} it follows that
\begin{align*}
\tilde\nabla_X'&\Psi_*f_*(A-\lambda I)DY
=\Psi_*\bar\nabla_Xf_*(A-\lambda I)DY
+\alpha^\Psi(f_*X, f_*(A-\lambda I)DY)\\
=&\,\Psi_*f_*\nabla_X(A-\lambda I)DY
+\<AX,(A-\lambda I)DY\>\Psi_*N
-\<X, (A-\lambda I)DY\>w\\
=&\,\Psi_*f_*(\nabla_X(A-\lambda I)D)Y
+\Psi_*f_*(A-\lambda I)D\nabla_XY\\
&+\<AX,(A-\lambda I)DY\>\Psi_*N
-\,\<X,(A-\lambda I)DY\>w\\
=&\,\Psi_*f_*(\nabla_X(A-\lambda I)D)Y
+\<D\nabla_XY,\nabla\lambda\>\Psi\circ f
-g_*\bar D\pi_*\nabla_XY\\
&+\,\<AX,(A-\lambda I)DY\>\Psi_*N
-\<X,(A-\lambda I)DY\>w.
\end{align*}
We obtain from the last two equations and
$\pi_*[X,Y]=[\bar{X},\bar{Y}]$ that
\begin{align*}
g_*((\nabla'_{\bar{Y}}\bar{D})\bar{X}&
-(\nabla'_{\bar{X}}\bar{D})\bar{Y})
+\a^g(\bar{Y},\bar{D}\bar{X})
-\a^g(\bar{X},\bar{D}\bar{Y})\\
&=\Psi_*f_*\Omega(X,Y)-\lambda\psi(X,Y)\Psi_*N
+\varphi(X,Y)\Psi\circ f+\psi(X,Y)w\nonumber
\end{align*}
where
\begin{align}
\Omega(X,Y)&=(\nabla_X(A-\lambda I)D)Y
-(\nabla_Y(A-\lambda I)D)X
-X\wedge Y(D^t\nabla\lambda)\label{om},\\
\psi(X,Y)&=\<Y,(A-\lambda I)DX\>
-\<X,(A-\lambda I)DY\>\label{ps},\\
\varphi(X,Y)&=\<(\nabla_YD)X
-(\nabla_XD)Y,\nabla\lambda\>+\hess\lambda(DX,Y)
-\hess\lambda(X,DY)\nonumber\\
-&\lambda(\<(A-\lambda I)X, (A-\lambda I)DY\>
-\<(A-\lambda I)DX,(A-\lambda I)Y\>).\label{vp}
\end{align}

It follows from \eqref{cod} and \eqref{gradients} 
that $\Omega$ vanishes. The symmetry of $\B$ yields 
$\psi=0$. Hence
$$
g_*((\nabla'_{\bar{Y}}\bar{D})\bar{X}
-(\nabla'_{\bar{X}}\bar{D})\bar{Y})
+\a^g(\bar{Y},\bar{D}\bar{X})
-\a^g(\bar{X},\bar{D}\bar{Y})
=\varphi(X,Y)\Psi\circ f.
$$
Since the term on the left-hand side is constant 
along the leaves of $\Delta$ then $\varphi$ has to 
vanish, which proves $\hyperlink{iv}{(iv)}$. 
\medskip

We prove the converse. Let $D=\mu J\in\Gamma(\End(\Delta^\perp))$ 
verify conditions $\hyperlink{i}{(i)}$ to $\hyperlink{v}{(v)}$.
In the sequel, we extend $D$ to an element of $\End(TM)$ 
defining $DT=0$ for any $T\in\Gamma(\Delta)$. Then 
$\hyperlink{v}{(v)}$ holds for any $X,Y\in\mathfrak{X}(M)$.

Set $F=\Psi\circ f\colon M^n\to\mathbb{V}^{n+2}\subset\Les^{n+3}$. 
Then let $\beta\colon TM\times TM\to N_F M$ be the 
symmetric tensor defined by
\be\label{defbeta}
\beta(X,Y)=\<(A-\lambda I)DX,Y\>(\Psi_*N+\lambda F)
\ee
where $N$ is a Gauss map of $f$. 
Let $B_\eta\in\Gamma(\End(TM))$ be given by 
$$
\<B_\eta X,Y\>=\<\beta(X,Y),\eta\>
$$ 
for any $\eta\in\Gamma(N_FM)$. For simplicity we write 
$N=\Psi_*N$. Observe that $B_N=(A-\lambda I)D$ and 
$B_w=\lambda B_N$. Since 
\be\label{sff2}
\a^F(X,Y)=\<AX,Y\>N-\<X,Y\>w,
\ee 
we have from $\hyperlink{v}{(v)}$ and  $A|_{\Delta}=\lambda I$ that
\be\label{derGauss}
A^F_{\beta(Y,Z)}X+B_{\a^F(Y,Z)}X-A^F_{\beta(X,Z)}Y-B_{\a^F(X,Z)}Y=0
\ee
for any $X,Y,Z\in\mathfrak{X}(M)$, where $A^F_\eta$ is the shape 
operator of $F$ with respect to $\eta\in\Gamma(N_FM)$.

We define $\mathcal{E}\colon TM\times N_FM\to N_FM$ by
\be\label{defE}
\mathcal{E}(X,N)=\<DX,\nabla \lambda\>F,\;\;
\mathcal{E}(X,w)=-\<DX,\nabla\lambda\>N\;\;\mbox{and}\;\;
\mathcal{E}(X,F)=0
\ee
for any $X\in\mathfrak{X}(M)$. Observe that $\mathcal{E}$ satisfies
the condition
\be\label{anti}
\<\mathcal{E}(X,\eta),\xi\>=-\<\mathcal{E}(X,\xi),\eta\>
\ee
for any $X\in\mathfrak{X}(M)$ and $\eta,\xi\in\Gamma(N_FM)$.

It follows from $\hyperlink{iii}{(iii)}$ that
\be\label{ones}
(\nabla_X B_N)Y-(\nabla_Y B_N)X
=\<DY,\nabla\lambda\>X-\<DX,\nabla\lambda\>Y
\ee
for any $X,Y\in\Gamma(\Delta^\perp)$.
Using $\hyperlink{ii}{(ii)}$ and $[D,C_T]=0$ we obtain
\begin{align*}
\!(\nabla_X B_N)T-(\nabla_T B_N)X&=B_NC_TX
\!-\!(\nabla_T(A-\lambda I))DX-(A-\lambda I)(\nabla_TD)X\\
&=(A-\lambda I)C_TDX-(\nabla_T(A-\lambda I))DX 
\end{align*}
for any $T\in\Gamma(\Delta)$. Now using the Codazzi equation, 
we have
\begin{align}\label{auxcod}
(\nabla_X B_N)T&-(\nabla_T B_N)X
= (A-\lambda I)C_TDX-(\nabla_{DX}A)T\nonumber\\
&=(A-\lambda I)C_TDX-\< DX,\nabla\lambda\>T
-(A-\lambda I)C_TDX\nonumber\\
&=-\< DX,\nabla\lambda\>T.
\end{align}
Since $\Delta$ is integrable, we obtain
\be\label{twos}
(\nabla_T B_N)S-(\nabla_S B_N)T=0
\ee
for any $T,S\in\Gamma(\Delta)$.
It follows from \eqref{ones}, \eqref{auxcod} and 
\eqref{twos} that
\be\label{derCod1}
(\nabla_X B_N)Y-(\nabla_Y B_N)X
=A^F_{\mathcal{E}(X,N)}Y-A^F_{\mathcal{E}(Y,N)}X
\ee
for any $X,Y\in\mathfrak{X}(M)$.

We have from \eqref{ones} that
\begin{align*}
(\nabla_X B_w)Y-(\nabla_Y B_w)X
&=\<X,\nabla\lambda\>B_NY
-\<Y,\nabla\lambda\>B_NX\\
&+\lambda \<DY,\nabla\lambda\>X-\lambda\<DX,\nabla\lambda\>Y
\end{align*}
for any $X,Y\in\Gamma(\Delta^\perp)$.
Let $\sigma\in\Gamma(\Delta^\perp)$ be given by 
$\nabla\lambda=(A-\lambda I)\sigma$. Using $\hyperlink{v}{(v)}$ 
we obtain
\begin{align}\label{three}
(\nabla_X& B_w)Y-(\nabla_Y B_w)X
=\<B_NY,\sigma\>(A-\lambda I)X
-\<B_NX,\sigma\>(A-\lambda I)Y\nonumber\\
&\;+\lambda \<DY,\nabla\lambda\>X
-\lambda\<DX,\nabla\lambda\>Y\nonumber\\
&=\<DY,\nabla\lambda\>(A-\lambda I)X
-\<DX,\nabla\lambda\>Y
+\lambda\<DY,\nabla\lambda\>X
-\lambda\<DX,\nabla\lambda\>Y\nonumber\\
&=\<DY,\nabla\lambda\>AX-\<DX,\nabla\lambda\>AY.
\end{align}
Using  \eqref{auxcod} it follows that
\begin{align}\label{four}
(\nabla_X B_w)T-(\nabla_T B_w)X
&=\lambda((\nabla_X B_N)T-(\nabla_T B_N)X)\nonumber\\
&=-\<DX,\nabla\lambda\>AT
\end{align}
for any $T\in\Gamma(\Delta)$.
As before, we have that
\be\label{five}
(\nabla_T B_w)S-(\nabla_S B_w)T=0
\ee
for any $S,T\in\Gamma(\Delta)$. We conclude from 
\eqref{three}, \eqref{four} and \eqref{five} that
\be\label{derCod2}
(\nabla_X B_w)Y-(\nabla_Y B_w)X
=A^F_{\mathcal{E}(X,w)}Y-A^F_{\mathcal{E}(Y,w)}X
\ee
for any $X,Y\in\mathfrak{X}(M)$.

We have that $B_F=0=\mathcal{E}(X,F)$, and hence
it holds trivially that
\be\label{derCod3}
(\nabla_X B_F)Y-(\nabla_Y B_F)X
=A^F_{\mathcal{E}(X,F)}Y-A^F_{\mathcal{E}(Y,F)}X
\ee
for any $X,Y\in\mathfrak{X}(M)$.

Next we focus on the covariant derivative of $\mathcal{E}$. 
Let $\nabla'^\perp$ denote the normal connection on $N_FM$.
We have
\begin{align*}
(\nabla'^\perp _X\mathcal{E})(Y,N)&
=\nabla'^\perp_X\mathcal{E}(Y,N)- \mathcal{E}(\nabla_XY,N)\\
&=X\<DY,\nabla\lambda\>F-\<D\nabla_XY,\nabla\lambda\>F\\
&=(\<(\nabla_XD)Y,\nabla\lambda\>+\hess\lambda(DY,X))F
\end{align*}
for any $X,Y\in\mathfrak{X}(M)$. Then
\begin{align*}
(\nabla'^\perp_X\mathcal{E})(Y,N)
-(\nabla'^\perp_Y\mathcal{E})(X,N)
&=(\<(\nabla_XD)Y-(\nabla_YD)X,\nabla\lambda\>\\
&+\hess\lambda(DY,X)-\hess\lambda(DX,Y))F
\end{align*}
for all $X,Y\in\mathfrak{X}(M)$.
From $\hyperlink{iv}{(iv)}$ we have
\be\label{six}
(\nabla'^\perp_X\mathcal{E})(Y,N)
-(\nabla'^\perp_Y\mathcal{E})(X,N)
=\lambda(\<B_NX,AY\>-\<AX,B_NY\>)F
\ee
for all $X,Y\in\Gamma(\Delta^\perp)$. 
Using $\hyperlink{ii}{(ii)}$ and $[D,C_T]=0$ we obtain
\begin{align}\label{seven}
(\nabla'^\perp_X\mathcal{E})(T,N)
&-(\nabla'^\perp_T\mathcal{E})(X,N)
=\mathcal{E}([T,X],N)-\nabla'^\perp_T\mathcal{E}(X,N)\nonumber\\
&=(\<DC_TX-(\nabla_TD)X,\nabla\lambda\>
-\hess\lambda(DX,T))F\nonumber\\
&=(\<C_TDX,\nabla\lambda\>
-\hess\lambda(DX,T))F\nonumber\\
&=(\<T,\nabla_{DX}\nabla\lambda\>
-\hess\lambda(DX,T))F\nonumber\\
&=0
\end{align}
for all $X\in\Gamma(\Delta^\perp)$ and $T\in\Gamma(\Delta)$.
We also have
\be\label{eigth}
(\nabla'^\perp_T\mathcal{E})(S,N)
-(\nabla'^\perp_S\mathcal{E})(T,N)=0
\ee
for any $S,T\in\Gamma(\Delta)$.

On the other hand, from \eqref{defbeta} and \eqref{sff2} 
we obtain
\begin{align}\label{nine}
\beta(X,AY)-\beta(AX,Y)&+\a^F(X,B_NY)-\a^F(B_NX,Y)\nonumber\\
&=\lambda(\<B_NX,AY\>-\<AX,B_NY\>)F
\end{align}
for any $X,Y\in\mathfrak{X}(M)$. From \eqref{six}, 
\eqref{seven}, \eqref{eigth} and \eqref{nine} we 
conclude that
\begin{align}\label{derRic1}
(\nabla'^\perp_X\mathcal{E})&(Y,N)
-(\nabla'^\perp_Y\mathcal{E})(X,N)\nonumber\\
&=\beta(X,AY)-\beta(AX,Y)+\a^F(X,B_NY)-\a^F(B_NX,Y)
\end{align}
for any $X,Y\in\mathfrak{X}(M)$.

Similarly as above, we obtain 
\begin{align*}
(\nabla'^\perp_X\mathcal{E})(Y,w)
-&(\nabla'^\perp_Y\mathcal{E})(X,w)
= \<(\nabla_YD)X-(\nabla_XD)Y,\nabla\lambda\>N\\
&+(\hess\lambda(DX,Y)-\hess\lambda(DY,X))N
\end{align*}
for all $X,Y\in\mathfrak{X}(M)$.
From $\hyperlink{iv}{(iv)}$ it follows that
$$
(\nabla'^\perp_X\mathcal{E})(Y,w)
-(\nabla'^\perp_Y\mathcal{E})(X,w)=
\lambda(\<AX,(A-\lambda I)DY\>-\<(A-\lambda I)DX,AY\>)N
$$
for $X,Y\in\Gamma(\Delta^\perp)$.
As before, we have from $\hyperlink{ii}{(ii)}$ and 
$[D,C_T]=0$ that
\begin{align*}
(\nabla'^\perp_X\mathcal{E})(T,w)
-(\nabla'^\perp_T\mathcal{E})(X,w)&=
(-\<C_TDX,\nabla\lambda\>+\hess\lambda(DX,T))N\\
&=(-\<T,\nabla_{DX}\lambda\>+\hess\lambda(DX,T))N\\
&=0
\end{align*}
and 
$$
(\nabla'^\perp_T\mathcal{E})(S,w)-(\nabla'^\perp_S\mathcal{E})(T,w)=0
$$
for any $T,S\in\Gamma(\Delta)$. It holds that
\begin{align*}
\beta(X,A^F_wY)-\beta(A^F_wX,Y)&+\a^F(X,B_wY)-\a^F(B_wX,Y)\\
&=\lambda(\<AX,B_NY\>-\<B_NX,AY\>)N
\end{align*}
for all $X,Y\in\mathfrak{X}(M)$. Thus
\begin{align}\label{derRic2}
(\nabla'^\perp_X\mathcal{E})&(Y,w)
-(\nabla'^\perp_Y\mathcal{E})(X,w)\nonumber\\
&=\beta(X,A^F_wY)-\beta(A^F_wX,Y)+\a^F(X,B_wY)-\a^F(B_wX,Y)
\end{align}
for all $X,Y\in\mathfrak{X}(M)$.
Finally, we have
$$
\beta(X,A^F_FY)-\beta(A^F_FX,Y)+\a^F(X,B_FY)-\a^F(B_FX,Y)=0
$$
for all $X,Y\in\mathfrak{X}(M)$, and since 
$\mathcal{E}(X,F)=0$, then
\begin{align}\label{derRic3}
(\nabla'^\perp_X\mathcal{E})&(Y,F)
-(\nabla'^\perp_Y\mathcal{E})(X,F)\nonumber\\
&=\beta(X,A^F_FY)-\beta(A^F_FX,Y)+\a^F(X,B_FY)-\a^F(B_FX,Y)
\end{align}
for all $X,Y\in\mathfrak{X}(M)$ holds trivially.

We have that $\beta$ is symmetric and  the tensor 
$\mathcal{E}$ satisfies condition \eqref{anti}. Moreover, 
the pair ($\mathcal{E},\beta)$ also satisfies \eqref{derGauss}, 
\eqref{derCod1}, 
\eqref{derCod2}, \eqref{derCod3}, \eqref{derRic1}, \eqref{derRic2} 
and \eqref{derRic3}. In this situation, the Fundamental theorem 
for infinitesimal variations, namely, Theorem $6$ 
in \cite{DJ1},  applies.  Notice that in the introduction of
\cite{DJ1} it was observed that Theorem $6$ holds for ambient
spaces of any signature, in particular, for the Lorentzian space
considered here.
Making also use of Proposition $5$ of \cite{DJ1}, we conclude 
that there is an infinitesimal bending 
$\tilde{\T}\in\Gamma(F^*(T\Les^{n+3}))$ of $F$ whose associated 
pair of tensors $(\tilde{\beta},\tilde{\mathcal{E}})$ satisfies
\be\label{tensors}
\tilde{\beta}=\beta+C\a^F\;\;\mbox{and}\;\;
\tilde{\mathcal{E}}=\mathcal{E}-\nabla^\perp C
\ee
where $C\in\Gamma(\End (N_FM))$ is skew-symmetric. 
Moreover, we have that $\tilde{\T}$ is unique up to 
trivial infinitesimal bendings. 

Write $\tilde{\T}$ as
$$
\tilde{\T}=\Psi_*\T+\<\tilde{\T},w\>F+\<\tilde{\T},F\>w.
$$
Being $\tilde{\T}$ an infinitesimal bending of $F$, 
we have
$$
\<\nab'_X\tilde{\T},F_*Y\>+\<\nab'_Y\tilde{\T},F_*X\>=0
$$
for all $X,Y\in\mathfrak{X}(M)$. Then
$$
\<\nab_X\T,f_*Y\>+\<\nab_Y\T,f_*X\>+2\<\tilde{\T},w\>\<X,Y\>=0
$$
for all $X,Y\in\mathfrak{X}(M)$.
Hence, setting $\rho=-\<\tilde{\T},w\>$ we have that 
$\T$ is a conformal infinitesimal bending of $f$ with 
conformal factor $\rho$.
It follows from \eqref{tensors} that the symmetric tensor 
$\B\in\Gamma(\End(TM))$ associated to $\T$ has 
the form $\B=B_N+cI$ where $c=-\<Cw,N\>$. And 
given that $B_N|_{\Delta^\perp}\neq 0$, we conclude that 
$\T$ is not trivial.

Any other conformal infinitesimal bending $\T'$ arising 
in this manner has the associated tensor $\B'=B_N+c'I$. 
Proposition \ref{trivial} now gives that $\T-\T'$ is trivial, 
and this concludes the proof.\qed

\begin{proposition}\label{parabolic}
Any parabolic  hypersurface $f\colon M^n\to \R^{n+1}$, 
$n\geq 5$, that admits a nontrivial conformal infinitesimal 
variation is conformally ruled.

Conversely, let $f\colon M^n\to\R^{n+1}$, $n\geq 5$, be 
a simply connected conformally ruled hypersurface  
free of points with a principal curvature of multiplicity 
at least $n-1$ and that is not conformally surface-like 
on any open subset of $M^n$. Then $f$ is parabolic
and admits a family of conformal infinitesimal 
bendings that are in one-to-one correspondence 
with the set of smooth functions on an interval.
\end{proposition}

\proof We have that $D=\mu J$ where $J^2=0$.  
Let $Y\in\Gamma(\Delta^\perp) $ be of unit-length 
such that $JY=0$ and let $X\in\Gamma(\Delta^\perp)$ be 
orthogonal to $Y$ satisfying $JX=Y$. That $\n_T^h J=0$ for 
any $T\in\Gamma(\Delta)$ is equivalent to 
\be\label{rul1}
\n_T^hY=0=\n_T^hX
\ee 
for all $T\in\Gamma(\Delta)$. Hence, replacing $J$ by 
$\|X\|J$, one can assume that also $X$ is of unit-length.

For the sequel, we extend $J$ to $TM$ as being zero on 
$\Delta$. Recall that 
$$
\B-bI=(A-\lambda I)D=\mu(A-\lambda I)J
$$ 
is symmetric. Then
\be\label{sym}
\<(A-\lambda I)Y,Y\>=\<(A-\lambda I)JX,Y\>=0.
\ee
Hence $(A-\lambda I)Y=\nu X$ where $\nu=\<AX,Y\>\neq 0$
by assumption. Then
$$
(\nabla_X\mu(A-\lambda I)J)Y-(\nabla_Y\mu(A-\lambda I)J)X
=-\mu(A-\lambda I)J\nabla_XY-\nabla_Y(\mu\nu X).
$$
On the other hand, we obtain from $\hyperlink{iii}{(iii)}$ 
that
\be\label{tres}
(\nabla_X\mu(A-\lambda I)J)Y
-(\nabla_Y\mu(A-\lambda I)J)X=-\mu Y(\lambda)Y.
\ee
Hence
$$
\mu(A-\lambda I)J\nabla_XY+\nabla_Y(\mu\nu X)
=\mu Y(\lambda)Y.
$$
Taking the inner product with $X$ and $Y$, 
respectively, gives
\be\label{aux}
Y(\mu\nu)=\mu\nu\<\nabla_XX,Y\>
\ee
and
\be\label{rul2}
Y(\lambda)=-\nu\<\nabla_YY,X\>.
\ee

Since $C_T\in\spa\{I,J\}$, we have 
\be\label{ct}
\<\nabla_YT,X\>=-\<C_TY,X\>=0
\ee
for any $T\in\Gamma(\Delta)$.
Let $T\in\Gamma(\Delta)$ be of unit length. The inner 
product of the Codazzi equation 
$(\nabla_TA)Y-(\nabla_YA)T=0$ with $T$ easily gives
\be\label{rul3}
Y(\lambda)=-\nu\<\nabla_TT,X\>.
\ee
It follows from \eqref{rul1}, \eqref{rul2}, \eqref{ct} and 
\eqref{rul3} that the subspaces  $\Delta\oplus\spa\{Y\}$ 
form an umbilical distribution.
Moreover,  we have from \eqref{sym} that 
$f$ restricted to any leaf of $\Delta\oplus\spa\{Y\}$ 
is umbilical in $\R^{n+1}$. Thus $f$ is conformally ruled.
\vspace{1ex}

We now prove the converse. Let $L$ be an $(n-1)$-dimensional 
umbilical distribution of $M^n$ such that the restriction of 
$f$ to any leaf is also umbilical. Therefore, there is 
$\lambda\in C^\infty(M)$ 
such that $L\subset\ker((A-\lambda I)_L)$, that is, 
$(A-\lambda I)(L)\subset L^\perp$. By assumption, we 
have that $\Delta=\ker(A-\lambda I)$ satisfies
$\dim\Delta=n-2$. 

Let  $X,Y$ be an orthonormal frame of $\Delta^\perp$ 
with $X$ orthogonal to $L$. Hence
\be\label{reg}
\<(A-\lambda I)Y,Y\>=0.
\ee
We have that $J\in\Gamma(\End(\Delta^\perp))$ defined 
by $JX=Y$  and $JY=0$  verifies $J^2=0$. It follows 
from \eqref{reg} that $(A-\lambda I)J$ is symmetric.
Now, since $L$ is umbilical, we have 
\be\label{par1}
\nabla_T^hY=0,
\ee
and this is equivalent to $\nabla_T^hJ=0$ for any 
$T\in\Gamma(\Delta)$. To verify that 
$C(\Gamma(\Delta))\subset\spa\{I,J\}$ it suffices to 
prove that $C_T\circ J=J\circ C_T$
for any $T\in\Gamma(\Delta)$. This is equivalent to
\be\label{equiva}
\<\nabla_YT,X\>=0\;\;\mbox{and}\;\;
\<\nabla_XX,T\>=\<\nabla_YY,T\>
\ee
for all $T\in\Gamma(\Delta)$. The first equation holds
since $L$ is umbilical. We have from \eqref{reg} that
\be\label{nu}
(A-\lambda I)Y=\nu X
\ee
where $\nu=\<AX,Y\>\neq 0$. From the Codazzi equation 
we easily obtain
$$
\nabla_T^hA=(A-\lambda I)C_T,
$$
and hence the right-hand side is symmetric.
We have
$$
\<(A-\lambda I)C_TX,Y\>=\nu\<\nabla_XX,T\>
\;\;\mbox{and}\;\;
\<(A-\lambda I)C_TY,X\>=\nu\<\nabla_YY,T\>,
$$
from where we obtain
$$
\<\nabla_XX,T\>=\<\nabla_YY,T\>
$$
for any $T\in\Gamma(\Delta)$.
Thus $f$ is parabolic with respect to $J$.

To show that $f$ admits a nontrivial conformal 
infinitesimal bending it suffices to prove that there 
is a smooth function $\mu$ such that the tensor 
$D=\mu J\in\Gamma(\End(\Delta^\perp))$ satisfies all 
conditions in Proposition \ref{fundD}. We already know that 
$(A-\lambda I)J$ is symmetric, hence condition $\hyperlink{i}{(i)}$ 
is satisfied for any function $\mu$. We assume that $\mu$ is 
constant along the leaves of $\Delta$, and now condition 
$\hyperlink{ii}{(ii)}$ follows from \eqref{par1}.
From the definition of $D$ it is easy to see that 
also condition $\hyperlink{v}{(v)}$ holds.

Condition $\hyperlink{iii}{(iii)}$ is just \eqref{tres}. 
We know that \eqref{rul3} holds for any 
$T\in\Gamma(\Delta)$ of unit-length. 
Hence and given that $L=\Delta\oplus\spa\{Y\}$ is an 
umbilical distribution, we obtain that \eqref{rul2} holds.  
But \eqref{rul2} is just the $Y$-component of \eqref{tres}.
The $X$-component of \eqref{tres} is \eqref{aux}, 
which can be stated as
$$
Y(\log\mu\nu)=\<\nabla_XX,Y\>.
$$
Choosing an arbitrary function as initial condition 
along one maximal integral curve of $X$, there exists 
a unique function $\mu$ such that $T(\mu)=0$ for all 
$T\in\Gamma(\Delta)$ and $\mu\nu$ is a solution of 
the preceding equation. Therefore, we have as many 
tensors $D$ satisfying $\hyperlink{iii}{(iii)}$ as 
smooth functions on an open interval.

We have that
$$
\<(\nabla_Y\mu J)X-(\nabla_X \mu J)Y,\nabla\lambda\>
=(Y(\mu)-\mu\<\nabla_XX,Y\>)Y(\lambda)
+\mu\<\nabla_YY,X\>X(\lambda).
$$
Choose any $D$ satisfying condition $\hyperlink{iii}{(iii)}$. 
Then \eqref{aux} and \eqref{rul2} yield
$$
\<(\nabla_Y\mu J)X-(\nabla_X \mu J)Y,\nabla\lambda\>
=-\frac{\mu}{\nu}Y(\lambda)(Y(\nu)+X(\lambda)).
$$
We have using  \eqref{rul2} that
\begin{align*}
\hess\lambda(\mu JX,Y)-\hess \lambda(X,\mu JY)
&=\mu(YY(\lambda)-\<\nabla_YY,X\>X(\lambda))\\
&=\mu(YY(\lambda)+\frac{1}{\nu}Y(\lambda)X(\lambda))
\end{align*}
and using \eqref{nu} that
$$
\lambda(\<(A-\lambda I)\mu JX,AY\>
-\<AX,(A-\lambda I)\mu JY\>)=\lambda\mu\nu^2.
$$
The last three equations give that the condition 
$\hyperlink{iv}{(iv)}$ is equivalent to
$$
YY(Y)-\frac{1}{\nu}Y(\lambda)Y(\nu)=-\lambda\nu^2,
$$
that can also be written as
\be\label{above}
Y((1/\nu)Y(\lambda))=-\lambda\nu .
\ee
To conclude, we show that \eqref{above} is just the 
Gauss equation
$$
\<R(Y,T)T,X\>=\<AT,T\>\<AY,X\>-\<AY,T\>\<AT,X\>
=\lambda\nu.
$$
In fact, we have using \eqref{rul3} and  \eqref{equiva} that
\begin{align*}
\<\nabla_Y\nabla_TT,X\>
&=Y\<\nabla_TT,X\>+\<\nabla_TT,Y\>\<\nabla_YY,X\>\\
&=-Y((1/\nu)Y(\lambda))+\<\nabla_TT,Y\>\<\nabla_YY,X\>.
\end{align*}
Also
$$
\<\nabla_T\nabla_YT,X\>=-\<\nabla_YT,\nabla_TX\>=0.
$$
Using \eqref{par1} we obtain
$$
\<\nabla_{[Y,T]}T,X\>=-\<\nabla_{\nabla_TY}T,X\>
=\<\nabla_TT,Y\>\<\nabla_TT,X\>.
$$
The last three equations yield 
$$
\<R(Y,T)T,X\>=-Y((1/\nu)Y(\lambda)).
$$
Now the proof follows from Proposition \ref{fundD}.\qed

\begin{proposition}\label{ruled}
Let $f\colon M^n\to\R^{n+1}$, $n\geq 5$, be a simply connected
conformally ruled hypersurface free of points with a
principal curvature of multiplicity at least $n-1$ and that is
not conformally surface-like on any open subset of $M^n$. Then
any conformal infinitesimal bending of $f$ is the variational
vector field of a conformal variation.
\end{proposition}

\proof We have seen that the conformal infinitesimal bendings of $f$
are in one-to-one correspondence with the tensors $D$ given  
in the proof of Proposition~ \ref{parabolic}. Take such a $D$ and let
$F\colon M^n\to\mathbb{V}^{n+2}\subset\Les^{n+3}$ be the
immersion $F=\Psi\circ f$, where $\Psi$ was given in \eqref{psi}. 
Let $\beta\colon TM\times TM\to N_FM$
and $\mathcal{E}\colon TM\times N_FM\to N_FM$ be given
by \eqref{defbeta} and \eqref{defE}, respectively. The
tensors $\beta$ and $\mathcal{E}$ are associated to an
infinitesimal bending of $F$, say $\Tilde{\T}$,
which determines a conformal infinitesimal bending $\T$ of $f$.
Let $\a^t\colon TM\times TM\to N_FM$, $t\in (-\epsilon,\epsilon)$,
be the symmetric tensor defined by
$$
\a^t(X,Y)=\a^F(X,Y)+t\beta(X,Y)
$$
for any $X,Y\in\mathfrak{X}(M)$. Since $\mathcal{E}$ satisfies
\eqref{anti} then
$\bar{\nabla}^{t}_X\eta=\nabla'^\perp_X\eta+t\mathcal{E}(X,\eta)$
is a connection on $N_FM$ that is compatible with the induced
metric, where $X\in\mathfrak{X}(M)$, $\eta\in\Gamma(N_FM)$ and
$\nabla'^\perp$ denotes the normal connection of $F$.

It follows from \eqref{derGauss}, \eqref{derCod1}, \eqref{derCod2},
\eqref{derCod3}, \eqref{derRic1}, \eqref{derRic2}, \eqref{derRic3}
together with the Gauss, Codazzi  and Ricci equations for $F$ that
$\a^t$ and $\bar{\nabla}^{t}$ verify the Gauss, Codazzi
and Ricci equations. Therefore, there is a family of isometric
immersions $F_t\colon M^n\to\mathbb{L}^{n+3}$ with $F_0=F$
together with vector bundle isometries
$\Phi_t\colon N_FM\to N_{F_t}M$ satisfying
$$
\a^{F_t}=\Phi_t\a^t\;\;\mbox{and}\;\;\nabla^{t\perp}\Phi_t
=\Phi_t(\bar{\nabla}^t)
$$
where $\a^{F_t}$ and $\nabla^{t\perp}$ are the second fundamental
form and normal connection of $F_t$, respectively. Then, we
have
$$
A^t_{\Phi_t F}X=-X\;\;\mbox{and}\;\;
\nabla_X^{t\perp}\Phi_t F=\Phi_t(\bar{\nabla}^t_XF)=0
$$
where $A^t_\eta$ is the shape operator of $F_t$ in the
direction of $\eta\in\Gamma(N_{F_t}M)$. Hence
$F_t-\Phi_t F=v_t$ is a constant vector field along $F_t$
for any $t$. Given that $\<F_t-v_t,F_t-v_t\>=0$, we obtain that
$F_t-v_t$ determines an isometric variation of $F_0=F$ in
$\mathbb{V}^{n+2}\subset\Les^{n+3}$. Hence, we assume that 
$F_t(x)\in\mathbb{V}^{n+2}$ for all $x\in M^n$.
The variational vector field 
$\tilde{\T}'=\partial/\partial_t|_{t=0}F_t$ is clearly an 
infinitesimal bending of $F$ and the tensor 
$\beta'$ associated to $\tilde{\T}'$ satisfies
$$
\beta'=(\partial/\partial_t|_{t=0}\a^{F_t})_{N_FM}
$$
(see the proof of Proposition 7 in \cite{DJ3}).
Since $\a^{F_t}=\Phi_t(\a+t\beta)$ then 
$$
\beta'=\beta+\Phi'\a^{F}
$$
where $\Phi'=\partial/\partial_t|_{t=0}\Phi_t\in\Gamma(\End(N_FM))$ 
is skew symmetric.

Let $\Pi\colon\mathbb{V}^{n+2}\setminus\R w\to\mathbb{E}^{m+1}
=\Psi(\R^{n+1})$ be the map $\Pi(u)=(1/\<u,v\>)u$. Then 
each $F_t$ induces an immersion $f_t\colon M^n\to\R^{n+1}$ 
such that $\Psi\circ f_t=\Pi\circ F_t$ for any $t$.
Observe that the metric induced by $f_t$ satisfies
$$
\<f_{t*}X,f{_t*}Y\>(x)
=\<(\Pi\circ F_t)_*X,(\Pi\circ F_t){_*}Y\>(x)
=\<F_t(x),w\>^{-2}\<X,Y\>(x)
$$
at any $x\in M^n$. Hence, the variation $f_t$ determines 
a conformal variation of $f$ in $\R^{n+1}$. The variational 
vector field $\T'$ is a conformal infinitesimal bending of 
$f$ with associated tensor
$\B'=B_N-\<\Phi'w,N\>I$.
Hence $\T-\T'$ is trivial, and this concludes the proof.
\vspace{2ex}\qed

\noindent\emph{Proof of Theorem \ref{main}:}
If $f$ admits a nontrivial conformal infinitesimal 
variation and is not conformally surface-like, we have from 
Proposition \ref{fundD} and Proposition 
\ref{parabolic} that $f$ is either hyperbolic or elliptic. 
The proof  of Proposition \ref{fundD} gives
that $D=\mu J$ is the lifting of a tensor 
$\bar{D}=\bar{\mu}\bar{J}$ on $L^2$.  Also
from that proof, we obtain 
\begin{align}\label{specials}
g_*((\nabla'_{\bar{Y}}\bar{D})\bar{X}&
-(\nabla'_{\bar{X}}\bar{D})\bar{Y})+\a^g(\bar{Y},\bar{D}\bar{X})
-\a^g(\bar{X},\bar{D}\bar{Y})\\
&=\Psi_*f_*\Omega(X,Y)-\lambda\psi(X,Y)\Psi_*N
+\varphi(X,Y)\Psi\circ f+\psi(X,Y)w \nonumber
\end{align}
where $X,Y\in\Gamma(\Delta^\perp)$ are the liftings 
of $\bar{X},\bar{Y}\in\mathfrak{X}(L^2)$ and $\Omega$, 
$\psi$ and $\varphi$ are given by \eqref{om}, \eqref{ps} 
and $\eqref{vp}$ respectively. Recall that $D$ satisfies 
the conditions $\hyperlink{i}{(i)}$ to $\hyperlink{v}{(v)}$.
Therefore, we have
\be\label{dcod}
(\nabla'_{\bar{X}}\bar{D})\bar{Y}=(\nabla'_{\bar{Y}}\bar{D})\bar{X}
\ee
and, since $\bar{D}=\bar{\mu}\bar{J}$, that
$$
\a^g(\bar{X},\bar{J}\bar{Y})=\a^g(\bar{J}\bar{X},\bar{Y}).
$$
Finally, that $g$ is a special hyperbolic or elliptic surface
follows from Proposition \ref{equivspecial} and the integrability 
condition of $\bar{\mu}$ in \eqref{dcod}.
\medskip

Conversely, let $f\colon M^n\to \R^{n+1}$ be 
parametrized by the conformal Gauss parametrization 
in terms of a special hyperbolic or a special elliptic 
pair. Then $f$ has a nowhere vanishing principal 
curvature $\lambda(x)$ at $x\in M^n$ of constant multiplicity 
$n-2$ and corresponding  eigenspace $\Delta(x)$.
Take $v=e_1$, $w=-2e_{n+3}$ and let 
$\Psi\colon \R^{n+1}\to\mathbb{V}^{n+2}\subset\Les^{n+3}$ 
be the embedding given by \eqref{psi}. Then  
$S\colon M^n\to\Sf_1^{n+2}$ given by \eqref{S} induces 
a map $s\colon L^2\to\Sf_1^{n+2}$ on the (local) space 
of leaves $L^2$ of $\Delta$. Moreover, by the choice of $v$ and 
$w$ we have that $s=g$.

We obtain from Proposition \ref{equivspecial} that, at least 
locally, there is a nowhere vanishing function 
$\bar{\mu}\in C^{\infty}(L^2)$ 
such that $\bar{D}=\bar{\mu}\bar{J}$ is a Codazzi tensor. 
Let $X,Y\in\Gamma(\Delta^\perp)$ be the liftings
of $\bar{X},\bar{Y}\in\mathfrak{X}(L)$.
If $D=\mu J$ is the lifting of $\bar{D}$ we have 
as before that \eqref{specials} holds. Given that $g$ 
is special hyperbolic or special elliptic, we 
have that 
$\Omega=\psi=\varphi=0$. In other words, we obtain that 
conditions $\hyperlink{i}{(i)}$, $\hyperlink{iii}{(iii)}$ 
and $\hyperlink{iv}{(iv)}$ are satisfied.

We recall that
$$
\<(\bar{D}\bar{X}\wedge\bar{Y}
-\bar{D}\bar{X}\wedge\bar{Y})\bar{Z},\bar{W}\>'=0
$$
for any $\bar{X},\bar{Y},\bar{Z},\bar{W}\in\mathfrak{X}(L)$. 
It follows from \eqref{metricsur} that
$$
\<((A-\lambda I)DX\wedge(A-\lambda I)Y
-(A-\lambda I)DY\wedge(A-\lambda I)X)
(A-\lambda I)Z,(A-\lambda I)W\>\!=\!0
$$
where $X,Y,Z,W\in\Gamma(\Delta^\perp)$ are the liftings 
of $\bar{X},\bar{Y},\bar{Z}$ and $\bar{W}$. Then
$$
(A-\lambda I)DX\wedge(A-\lambda I)Y
-(A-\lambda I)DY\wedge(A-\lambda I)X=0
$$
for all $X,Y\in\Gamma(\Delta^\perp)$, and hence 
$\hyperlink{v}{(v)}$ holds. Given that $D$ is projectable
it follows from Corollary $11.7$ of \cite{DT2} that
$\nabla_T^hD=[D,C_T]=0$ for all $T\in\Gamma(\Delta)$. Hence
$\hyperlink{ii}{(ii)}$ holds. Now the proof follows 
from Proposition \ref{fundD}.\qed

\begin{remark}\label{cartan} 
{\em In order to obtain, in terms of the conformal Gauss 
parametrization, that a nontrivial  conformal infinitesimal bending 
is, in fact, the variational vector field of a conformal variation 
one has to require the special hyperbolic or special elliptic surface 
to satisfy a strong additional condition, namely, 
that $\gg_u=\gh_v=2\gg\gh$ in the former case and 
$\Gamma_z = 2\Gamma\bar{\Gamma}$
in the latter case, see \cite{DT1} or \cite{DT2}.}
\end{remark}

\noindent\emph{Proof of Theorem \ref{main2}:} 
The proof follows from Propositions \ref{parabolic}
and \ref{ruled}.\vspace{2ex}\qed

Part of this work is the result of the visit 21171/IV/19 funded 
by the Fundación Séneca-Agencia de Ciencia y Tecnología de la
Región de Murcia in connection with the ``Jiménez De La Espada"
Regional Programme For Mobility, Collaboration And Knowledge Exchange.

Marcos Dajczer was partially supported by the Fundación Séneca Grant\\
21171/IV/19 (Programa Jiménez de la Espada), MICINN/FEDER project
PGC2018-097046-B-I00, and Fundación Séneca project 19901/GERM/15, 
Spain.

Miguel I. Jimenez thanks the mathematics department of the Universidad 
de Murcia for the hospitality during his visit where part of this 
work was developed.

Theodoros Vlachos acknowledges support 
by the General Secretariat for Research and Technology 
(GSRT) and the Hellenic Foundation for Research and 
Innovation (HFRI) Grant No: 133.

\noindent Marcos Dajczer\\
IMPA -- Estrada Dona Castorina, 110\\
22460--320, Rio de Janeiro -- Brazil\\
e-mail: marcos@impa.br

\bigskip

\noindent Miguel Ibieta Jimenez\\
Instituto de Ciências Matemáticas e de Computação\\
Universidade de São Paulo\\
São Carlos\\
SP 13560-970-- Brazil\\
e-mail: mibieta@impa.br

\bigskip

\noindent Theodoros Vlachos\\
University of Ioannina \\
Department of Mathematics\\
Ioannina -- Greece\\
e-mail: tvlachos@uoi.gr


\begin{thebibliography}{lll}

\bibitem{Ca1} Cartan, E.,
\emph{La d\'eformation des hypersurfaces dans 
l'espace euclidien r\'eel a $n$ dimensions},
Bull. Soc. Math. France {\bf 44} (1916), 65--99.

\bibitem{Ca2}  Cartan, E., 
\emph{La d\'eformation des
hypersurfaces dans l'espace conforme r\'eel a $n\geq 5$ 
dimensions},
Bull. Soc. Math. France {\bf 45} (1917), 57--121.

\bibitem{Ca3} Cartan, E., 
\emph{Sur certains hypersurfaces de l'espace conforme r\'eel a 
cinq dimensions}, 
Bull. Soc. Math. France 46 (1918), 84--105.

\bibitem{DF} Dajczer, M. and Florit, L., 
\emph{Compositions of isometric immersions in higher codimension},
Manuscripta Math. {\bf 105} (2001), 507--517. 

\bibitem{DJ1} Dajczer, M. and Jimenez, M. I., 
\emph{Infinitesimal variations of submanifolds}, 
to appear in Bull. Braz. Math. Soc. 
\url{https://arxiv.org/abs/1911.01863}

\bibitem{DJ2} Dajczer, M. and Jimenez, M. I., 
\emph{Conformal infinitesimal variations of submanifolds}, 
preprint. \url{https://arxiv.org/abs/2002.02551}

\bibitem{DJ3} Dajczer, M. and Jimenez, M. I.,
\emph{Genuine infinitesimal bendings of submanifolds},
preprint. \url{https://arxiv.org/abs/1904.10409}

\bibitem{DT1} Dajczer, M. and  Tojeiro, R.,
\emph{On Cartan's conformally deformable hypersurfaces},
Michigan Math. J. {\bf 47}  (2000), 529--557.

\bibitem{DT2} Dajczer, M. and  Tojeiro, R., 
``Submanifold theory beyond an introduction".
Universitext. Springer, 2019.

\bibitem{DV} Dajczer, M. and Vlachos, Th.,
\emph{The infinitesimally bendable Euclidean hypersurfaces},
Ann. Mat. Pura Appl. {\bf 196} (2017), 1961--1979 and 
Ann. Mat. Pura Appl. {\bf 196} (2017), 1981--1982.

\bibitem{Ji} Jimenez, M. I., 
\emph{Infinitesimal bendings of complete Euclidean 
hypersurfaces}, 
Manuscripta Math. {\bf 157} (2018), 513--527.

\bibitem{Sb1} Sbrana, U.,
\emph{Sulle variet\`a ad $n-1$ dimensioni deformabili nello 
spazio euclideo ad $n$ dimensioni},
Rend. Circ. Mat. Palermo  {\bf 27} (1909), 1--45.

\bibitem{Sb2} Sbrana, U.,
\emph{Sulla deformazione infinitesima delle ipersuperficie},
Ann. Mat. Pura  Appl. {\bf 15} (1908), 329--348.

\bibitem{Sc} Schottenloher, M., 
``A mathematical introduction to conformal field theory".
Lecture Notes in Physics, 759. Springer-Verlag, 2008.

\end{thebibliography}
\end{document}